\documentclass[reqno,centertags, 12pt]{amsart}

\usepackage{amssymb,amsmath,amsfonts,amssymb}
-\marginparwidth \oddsidemargin -4mm \evensidemargin \oddsidemargin
\usepackage{latexsym, verbatim}
\advance\hoffset by 5mm

\def\@abssec#1{\vspace{.05in}\footnotesize \parindent .2in
{\bf #1. }\ignorespaces}
\newtheorem{theorem}{Theorem}[section]
\newtheorem{thm}[theorem]{Theorem}
\newtheorem{lemma}[theorem]{Lemma}
\newtheorem{proposition}[theorem]{Proposition}
\newtheorem{corollary}[theorem]{Corollary}
\newtheorem{definition}[theorem]{Definition}

\def \Rm {\mathbb R}

\def \Zm {\mathbb Z}

\newcommand{\eps}{\varepsilon}

\newcommand{\del}{\delta}
\newcommand{\tht}{\theta}

\newcommand{\Tm}{{\mathbb T}}
\newcommand{\Sm}{{\mathbb S}}

\DeclareMathOperator{\Lip}{Lip}

\newcommand{\calP}{\mathcal{P}}
\allowdisplaybreaks \numberwithin{equation}{section}

\title{Diffusion and Mixing in Fluid Flow}
\author{Peter Constantin}
\thanks{Department of
Mathematics, University of Chicago, Chicago, IL 60637; email:
const@cs.uchicago.edu}
\author{Alexander Kiselev}
\thanks{Department of
Mathematics, University of Wisconsin, Madison, WI 53706; e-mail:
kiselev@math.wisc.edu}
\author{Lenya Ryzhik}
\thanks{Department of
Mathematics, University of Chicago, Chicago, IL 60637; email:
ryzhik@math.uchicago.edu}
\author{Andrej Zlato\v s}
\thanks{Department of
Mathematics, University of Wisconsin, Madison, WI 53706; e-mail:
zlatos@math.wisc.edu}

\begin{document}

\begin{abstract}
We study enhancement of diffusive mixing on a compact Riemannian
manifold by a fast incompressible flow. Our main result is a sharp
description of the class of flows that make the deviation of the
solution from its average arbitrarily small in an arbitrarily
short time, provided that the flow amplitude is large enough. The
necessary and sufficient condition on such flows is expressed
naturally in terms of the spectral properties of the dynamical
system associated with the flow. In particular, we find that
weakly mixing flows always enhance dissipation in this sense. The
proofs are based on a general criterion for the decay of the
semigroup generated by an operator of the form $\Gamma+iAL$ with a
negative unbounded self-adjoint operator $\Gamma$, a self-adjoint
operator $L$, and parameter $A\gg 1$. In particular,  they employ
the RAGE theorem describing evolution of a quantum state belonging
to the continuous spectral subspace of the hamiltonian (related to
a classical theorem of Wiener on Fourier transforms of measures).
Applications to quenching in reaction-diffusion equations are also considered.
\end{abstract}

\maketitle

\section{Introduction}

Let $M$ be a smooth compact $d$-dimensional Riemannian manifold.
The main objective of this paper is the study of the effect of a
strong incompressible flow on diffusion on $M.$ Namely, we
consider solutions of the passive scalar equation
\begin{equation}\label{maineq}
\phi^A_t(x,t) +A u \cdot \nabla \phi^A(x,t) - \Delta
\phi^A(x,t)=0, \,\,\,\,\,\phi^A(x,0)=\phi_0(x).
\end{equation}
Here $\Delta$ is the Laplace-Beltrami operator on $M,$ $u$ is a
divergence free vector field, $\nabla$ is the covariant
derivative,
and $A\in\Rm$ is a parameter regulating the strength of the flow.
We are interested in the behavior of solutions of (\ref{maineq})
for $A\gg 1$ at a fixed time $\tau>0$.

It is well known that as time tends to infinity, the solution
$\phi^A(x,t)$ will tend to its average,
\[ \overline{\phi} \equiv \frac{1}{|M|}\int\limits_M \phi^A(x,t)\,d\mu=
\frac{1}{|M|}\int\limits_M \phi_0(x)\,d\mu, \] with $|M|$ being the
volume of $M$.
We would like to understand how the speed of convergence to the
average depends on the properties of the flow and determine which
flows are efficient in enhancing the relaxation process.

The question of the influence of advection on diffusion is very
natural and physically relevant, and the subject has a long
history. The passive scalar model is one of the most studied
PDEs in both mathematical and physical literature. One important
direction of research focused on homogenization, where in a long
time--large propagation distance limit the solution of a passive
advection-diffusion equation converges to a solution of an effective
diffusion equation. Then one is interested in the dependence of the
diffusion coefficient on the strength of the fluid flow. We refer to
\cite{KM} for more details and references. The main difference with the
present work is that here we are interested in the flow effect in a finite time
without the long time limit.

On the other hand, the Freidlin-Wentzell theory
\cite{Freidlin1,Freidlin2,FW1,FW2} studies \eqref{maineq} in $\Rm^{2}$
and, for a class of Hamiltonian flows, proves the convergence of
solutions as $A \rightarrow \infty$ to solutions of an effective
diffusion equation on the Reeb graph of the hamiltonian. The graph, essentially, is obtained by
identifying all points on any streamline. The conditions on the flows
for which the procedure can be carried out are given in terms of
certain non-degeneracy and growth assumptions on the stream
function. The Freidlin-Wentzell method does not apply, in particular,
to ergodic flows or in odd dimensions.

Perhaps the closest to our setting is the work of Kifer and more
recently a result of Berestycki, Hamel and Nadirashvili. Kifer's work
(see \cite{Kifer1, Kifer2, Kifer3,Kifer4} where further references can
be found) employs probabilistic methods and is focused, in particular,
on the estimates of the principal eigenvalue (and, in some special
situations, other eigenvalues) of the operator $-\epsilon \Delta + u
\cdot \nabla$ when $\epsilon$ is small, mainly in the case of the
Dirichlet boundary conditions. In particular, the asymptotic behavior
of the principal eigenvalue $\lambda^\epsilon_0$ and the corresponding
positive eigenfunction $\phi^\epsilon_0$ for small $\epsilon$ has been
described in the case where the operator $u \cdot \nabla$ has a
discrete spectrum and sufficiently smooth eigenfunctions. It is well
known that the principal eigenvalue determines the asymptotic rate of
decay of the solutions of the initial value problem, namely
\begin{equation}\label{asbeh}
\lim_{t \rightarrow \infty} t^{-1}\log \|\phi^\epsilon
(x,t)\|_{L^2} = -\lambda^\epsilon_0
\end{equation}
(see e.g. \cite{Kifer2}). In a related recent work \cite{BHN},
Berestycki, Hamel and Nadirashvili utilize PDE methods to prove a
sharp result on the behavior of the principal eigenvalue
$\lambda_A$ of the operator $-\Delta + A u \cdot \nabla$ defined
on a bounded domain $\Omega \subset \Rm^d$ with the Dirichlet
boundary conditions. The main conclusion is that $\lambda_A$ stays
bounded as $A \rightarrow \infty$ if and only if $u$ has a first
integral $w$ in $H^1_0(\Omega)$ (that is, $u \cdot \nabla w =0$). An
elegant variational principle determining the limit of $\lambda_A$
as $A \rightarrow \infty$ is also proved. In addition, \cite{BHN}
provides a direct link between the behavior of the principal
eigenvalue and the dynamics which is more robust than
\eqref{asbeh}: it is shown that $\|\phi^A(\cdot,1)\|_{L^2(\Omega)}$ can
be made arbitrarily small for any initial datum by increasing $A$
if and only if $\lambda_A \rightarrow \infty$ as $A \rightarrow
\infty$ (and, therefore, if and only if the flow $u$ does not have
a first integral in $H^1_0(\Omega)$). We should mention that there
are many earlier works providing variational characterization of
the principal eigenvalues, and refer to \cite{BHN,Kifer4} for more
references.

Many of the studies mentioned above also apply in the case of a
compact manifold without boundary or Neumann boundary conditions,
which is the primary focus of this paper. However, in this case the
principal eigenvalue is simply zero and corresponds to the constant
eigenfunction. Instead one is interested in the speed of convergence
of the solution to its average, the relaxation speed. A recent work of
Franke \cite{Franke} provides estimates on the heat
kernels corresponding to the incompressible drift and diffusion on
manifolds, but these estimates lead to upper bounds on
$\|\phi^A(1)-\overline{\phi}\|$ which essentially do not improve as $A
\rightarrow \infty.$ One way to study the convergence speed is to
estimate the spectral gap -- the difference between the principal
eigenvalue and the real part of the next eigenvalue. To the best of
our knowledge, there is very little known about such estimates in the
context of \eqref{maineq}; see \cite{Kifer2} p.  251 for a
discussion. Neither probabilistic methods nor PDE methods of
\cite{BHN} seem to apply in this situation, in particular because the
eigenfunction corresponding to the eigenvalue(s) with the second
smallest real part is no longer positive and the eigenvalue itself
does not need to be real.  Moreover, even if the spectral gap estimate
were available, generally it only yields a limited asymptotic in time
dynamical information of type \eqref{asbeh}, and how fast the long time limit is
achieved may depend on $A.$ Part of our motivation for studying the
advection-enhanced diffusion comes from the applications to quenching
in reaction-diffusion equations (see e.g.
\cite{CKR,FKR,KZ,Roq,Z}), which we discuss in
Section~\ref{quenching}. For these applications, one needs
estimates on the $A$-dependent $L^\infty$ norm decay at a fixed
positive time, the type of information the bound like
\eqref{asbeh} does not provide. We are aware of only one case
where enhanced relaxation estimates of this kind are available. It
is the recent work of Fannjiang, Nonnemacher and
Wolowski \cite{FW,FNW}, where such estimates are provided in the
discrete setting (see also \cite{Kifer2} for some related earlier
references). In these papers a unitary evolution step
(a certain measure preserving map on the torus)
alternates with a dissipation step, which, for example, acts
simply by multiplying the Fourier coefficients by damping factors.
The absence of sufficiently regular eigenfunctions appears as a
key for the lack of enhanced relaxation in this particular class
of dynamical systems. In \cite{FW,FNW}, the authors also provide
finer estimates of the dissipation time for particular classes of
toral automorphisms (that is, they estimate how many steps are
needed to reduce the $L^2$ norm of the solution by a factor of two
if the diffusion strength is $\epsilon$).

Our main goal in this paper is to provide a sharp characterization
of incompressible flows that are relaxation enhancing, in a quite
general setup. We work directly with dynamical estimates, and do not
discuss the spectral gap. The following natural definition will be
used in this paper as a measure of the flow efficiency in improving
the solution relaxation.

\begin{definition}\label{relen}
Let $M$ be a smooth compact Riemannian manifold. The incompressible
flow $u$ on $M$ is called relaxation enhancing if for every $\tau>0$
and $\delta>0,$ there exist $A(\tau, \delta)$ such that for any
$A>A(\tau, \delta)$ and any $\phi_0\in L^2(M)$ with
$\|\phi_0\|_{L^2(M)}= 1$ we have
\begin{equation}\label{disscond}
\|\phi^A(\cdot,\tau)-\overline{\phi}\|_{L^2(M)} < \delta,
\end{equation}
where $\phi^A(x,t)$ is the solution of \eqref{maineq} and
$\overline\phi$ the average of $\phi_0$.
\end{definition}

\noindent {\it Remarks.}
1. In Theorem \ref{anyp} we show that the choice of the $L^2$ norm in the
definition is not essential and can be replaced by any $L^p$-norm with $1\le p\le\infty$.
\smallskip

2. It follows from the proofs of our main results that the
relaxation enhancing class is not changed even when we allow
the flow strength that ensures \eqref{disscond} to depend on
$\phi_0$, that is, if we require \eqref{disscond} to hold for all
$\phi_0\in L^2(M)$ with $\|\phi_0\|_{L^2(M)}= 1$ and all
$A>A(\tau,\delta,\phi_0)$.
\smallskip

Our first result is as follows.
%
%
\begin{thm}\label{fluid}
Let $M$ be a smooth compact Riemannian manifold. A Lipschitz continuous
incompressible flow $u \in \Lip(M)$ is relaxation enhancing if and
only if the operator $u \cdot \nabla$
has no eigenfunctions in $H^1(M),$ other than the constant
function.
\end{thm}

Any incompressible flow $u \in \Lip(M)$ generates a unitary
evolution group $U^t$ on $L^2(M),$ defined by $U^tf(x) =
f(\Phi_{-t}(x)).$ Here $\Phi_t(x)$ is a measure preserving
transformation associated with the flow, defined by $\frac{d}{dt}
\Phi_t(x) = u(\Phi_t(x)),$ $\Phi_0(x)=x.$ Recall that a flow $u$
is called weakly mixing if the corresponding operator $U$ has only
continuous spectrum. 
The weakly mixing flows are
ergodic, but not necessarily mixing (see e.g. \cite{CFS}). There
exist fairly explicit examples of weakly mixing flows
\cite{Anosov,Fayad1, Fayad2,Kolmogorov,Sklover,vN}, some of which
we will discuss in Section~\ref{flowex}. A direct consequence of
Theorem~\ref{fluid} is the following Corollary.

\begin{corollary}\label{wm}
Any weakly mixing incompressible flow $u \in \Lip(M)$ is relaxation
enhancing.
\end{corollary}

Theorem~\ref{fluid}, as we will see in Section~\ref{flows}, in its
turn follows from quite general abstract criterion, which we are now
going to describe. Let $\Gamma$ be a self-adjoint, positive,
unbounded operator with a discrete spectrum on a separable Hilbert
space $H.$ Let $0 < \lambda_1 \leq \lambda_2 \leq \dots$ be the
eigenvalues of $\Gamma,$ and $e_j$ the corresponding orthonormal
eigenvectors forming a basis in $H.$ The (homogenous) Sobolev space
$H^m(\Gamma)$ associated with $\Gamma$ is formed by all vectors
$\psi = \sum_j c_j e_j$ such that
\[
\|\psi\|_{H^m(\Gamma)}^2 \equiv \sum\limits_j \lambda_j^{m} |c_j|^2
< \infty.
\]
Note that $H^2(\Gamma)$ is the domain $D(\Gamma)$ of $\Gamma$. Let
$L$ be a self-adjoint operator such that, for any $\psi \in
H^1(\Gamma)$ and $t>0$ we have
\begin{equation}\label{dyncond}
\|L\psi\|_H \leq C \|\psi\|_{H^1(\Gamma)} \,\,\,{\rm and}\,\,\,
\|e^{iLt} \psi\|_{H^1(\Gamma)} \leq B(t) \|\psi\|_{H^1(\Gamma)}
\end{equation}
with both the constant $C$ and the function $B(t)<\infty$ independent
of $\psi$ and $B(t)\in L^2_{\rm loc}(0,\infty)$. Here $e^{iLt}$ is the
unitary evolution group generated by the self-adjoint operator $L.$
One might ask whether one of the two conditions in
\eqref{dyncond} does not imply the other. We show at the end of
Section \ref{prelim}, by means of an example, that this is not the
case in general.

Consider a solution $\phi^A(t)$ of the Bochner differential equation
\begin{equation}\label{hseq}
\frac{d}{dt}\phi^A(t) = iAL \phi^A(t) - \Gamma \phi^A(t),
\,\,\,\,\,\, \phi^A(0)=\phi_0.
\end{equation}
\begin{thm}\label{general}
Let $\Gamma$ be a self-adjoint, positive, unbounded operator with a
discrete spectrum and let a self-adjoint operator $L$ satisfy
conditions \eqref{dyncond}. Then the following two statements are
equivalent:
\begin{itemize}
\item For any $\tau,\delta >0$ there exists $A(\tau, \delta)$
such that for any $A> A(\tau, \delta)$ and any $\phi_0\in H$ with
$\|\phi_0\|_H=1$, the solution $\phi^A(t)$ of the equation
\eqref{hseq} satisfies $\|\phi^A(\tau)\|_H^2 < \delta.$ \item The
operator $L$ has no eigenvectors lying in $H^1(\Gamma).$
\end{itemize}
\end{thm}

\noindent {\it Remark.} Here $L$ corresponds to
$iu\cdot\nabla$ (or, to be precise, a self-adjoint operator
generating the unitary evolution group $U^t$ which is equal to
$iu\cdot\nabla$ on $H^1(M)$), and $\Gamma$ to $-\Delta$ in Theorem
\ref{fluid}, with $H\subset L^2(M)$ the subspace of mean zero
functions.
\smallskip

Theorem~\ref{general} provides a sharp answer to the general
question of when a combination of fast unitary evolution and
dissipation produces a significantly stronger dissipative effect
than dissipation alone. It can be useful in any model describing a
physical situation which involves fast unitary dynamics with
dissipation (or, equivalently, unitary dynamics with weak
dissipation). We prove Theorem~\ref{general} in Section~\ref{abs}.
The proof uses ideas from quantum dynamics, in particularly the
RAGE theorem (see e.g., \cite{CFKS}) describing evolution of a
quantum state belonging to the continuous spectral subspace of a
self-adjoint operator.

A natural concern is if the existence of rough eigenvectors of $L$
is consistent with the condition \eqref{dyncond} which says that
the dynamics corresponding to $L$ preserves $H^1(\Gamma)$. In
Section~\ref{rough} we answer this question in the affirmative by
providing examples where rough eigenfunctions exist yet
\eqref{dyncond} holds. One of them involves a discrete version of
the celebrated Wigner-von Neumann construction of an imbedded
eigenvalue of a Schr\"odinger operator \cite{WvN}. Moreover, in
Section~\ref{flowex} we describe an example of a smooth flow on
the two dimensional torus $\Tm^2$ with discrete spectrum and rough
(not $H^1(\Tm^2)$) eigenfunctions -- this example essentially goes
back to Kolmogorov \cite{Kolmogorov}. Thus, the result of
Theorem~\ref{general} is precise.

In  Section~\ref{quenching}, we discuss the application of
Theorem~\ref{fluid} to quenching for reaction-diffusion equations
on compact manifolds and domains. This corresponds to adding a
non-negative reaction term $f(T)$ on the right hand side of
\eqref{maineq}, with $f(0)=f(1)=0.$ Then the long-term dynamics
can lead to two outcomes: $\phi^A \rightarrow 1$ at every point
(complete combustion), or $\phi^A \rightarrow c <1$ (quenching).
The latter case is only possible if $f$ is of the ignition type,
that is, there exists $\theta_0$ such that $f(T) =0$ for $T
\leq \theta_0,$ and $c \leq \theta_0.$ The question is then how
the presence of strong fluid flow may aid the quenching process.
We note that quenching/front propagation in infinite domains is
also of considerable interest. Theorem~\ref{fluid} has
applications in that setting as well, but they will
be considered elsewhere.


\section{Preliminaries}\label{prelim}

In this section we collect some elementary facts and estimates for the
equation \eqref{hseq}.  Henceforth we are going to denote the standard
norm in the Hilbert space $H$ by $\|\cdot\|,$ the inner product in $H$
by $\langle \cdot,\cdot \rangle,$ the Sobolev spaces $H^m(\Gamma)$
simply by $H^m$ and norms in these Sobolev spaces by $\|\cdot\|_m.$ We
have the following existence and uniqueness theorem.

\begin{thm}\label{exun}
Assume that for any $\psi \in H^1,$ we have
\begin{equation}\label{Lcon1}
\|L\psi\| \leq C \|\psi\|_{1}.
\end{equation}
 Then for any
$T>0,$ there exists a unique solution $\phi(t)$ of the equation
\[
\phi'(t) = (iL - \Gamma)\phi(t), \,\,\,\,\phi(0)=\phi_0 \in H^1.
\]
This solution satisfies
\begin{equation}\label{reg}
\phi(t) \in L^2([0,T], H^2) \cap C([0,T],H^1), \,\,\, \phi'(t) \in
L^2([0,T], H).
\end{equation}
\end{thm}

{\it Remarks.} 1. The proof of Theorem~\ref{exun} is standard, and
can proceed by constructing a weak solution using Galerkin
approximations and then establishing uniqueness and regularity. We
refer, for example, to Evans \cite{Evans} where the construction is
carried out for parabolic PDEs but, given the assumption
\eqref{Lcon1}, can be applied verbatim in the general case.
\smallskip

2. The existence theorem is also valid for initial data $\phi_0 \in
H,$ but the solution has rougher properties at intervals containing
$t=0,$ namely
\begin{equation}\label{wreg}
\phi(t) \in L^2([0,T], H^1) \cap C([0,T],H), \,\,\, \phi'(t) \in
L^2([0,T], H^{-1}).
\end{equation}
The existence of a rougher solution can also be derived from the
general semigroup theory, by checking that $iL-\Gamma$ satisfies
the conditions of the Hille-Yosida theorem and thus generates a
strongly continuous contraction semigroup
in $H$ (see, e.g. \cite{EN}).
\smallskip

Next we establish a few properties that are more specific to our
particular problem. It will be more convenient for us, in terms of
notation, to work with an equivalent reformulation of
\eqref{hseq}, by setting $\epsilon = A^{-1}$ and rescaling time by
the factor $\epsilon^{-1}$, thus arriving at the equation
\begin{equation}\label{sdeq}
(\phi^\epsilon)'(t) = (iL-\epsilon \Gamma)\phi^\epsilon(t),
\,\,\,\phi^\epsilon(0) = \phi_0.
\end{equation}

\begin{lemma}\label{gradbound}
Assume \eqref{Lcon1}, then for any initial data $\phi_0 \in H,$ $\|\phi_0\|=1,$
the solution $\phi^\epsilon(t)$ of \eqref{sdeq} satisfies
\begin{equation}\label{gcont} \epsilon \int\limits_0^\infty
\|\phi^\epsilon(t)\|_1^2 dt \leq \frac12.
\end{equation}
\end{lemma}
\begin{proof}
Recall that if $\phi\in H^1(\Gamma)$, then $\Gamma\phi\in
H^{-1}(\Gamma)$ and $\langle \Gamma \phi,\phi\rangle =
\|\phi\|_{1}^2$. The regularity conditions
\eqref{reg}-\eqref{wreg} and the fact that $L$ is self-adjoint
allow us to compute
\begin{equation}\label{dissip}
\frac{d}{dt} \|\phi^\epsilon\|^2 = \langle \phi^\epsilon,
\phi^\epsilon_t \rangle + \langle \phi^\epsilon_t, \phi^\epsilon
\rangle =
-2\epsilon \|\phi^\epsilon \|^2_1. \end{equation} Integrating in
time and taking into account the normalization of $\phi_0,$ we
obtain (\ref{gcont}).
\end{proof}

An immediate consequence of \eqref{dissip} is the following
result, that we state here as a separate lemma for convenience.

\begin{lemma}\label{easy}
Suppose that for all times $t\in(a,b)$ we have
$\|\phi^\epsilon(t)\|_1^2 \geq N \|\phi^\epsilon(t)\|^2.$ Then the
following decay estimate holds:
\[ \|\phi^\epsilon(b)\|^2 \leq e^{-2\epsilon N (b-a)}
\|\phi^\epsilon(a)\|^2. \]
\end{lemma}

Next we need an estimate on the growth of the difference between
solutions corresponding to $\epsilon>0$ and $\epsilon=0$ in the
$H$-norm.

\begin{lemma}\label{approxH}
Assume, in addition to \eqref{Lcon1}, that for any $\psi \in H^1$
and $t>0$ we have
\begin{equation}\label{Lcon2}
\|e^{iLt}\psi\|_1 \leq B(t)\|\psi\|_1
\end{equation}
for some $B(t)<\infty$ such that $B(t)\in L^2_{\rm
loc}[0,\infty).$ Let $\phi^0(t),\phi^\epsilon(t)$ be solutions of
\[
(\phi^0)'(t) = iL\phi^0(t), \,\,\, (\phi^\epsilon)'(t) =
(iL-\epsilon \Gamma)\phi^\epsilon(t),
\]
satisfying $\phi^0(0)=\phi^\epsilon(0) = \phi_0 \in H^1.$ Then we
have
\begin{equation}\label{l2approx}
\frac d{dt} \|\phi^\epsilon (t) - \phi^0 (t)\|^2  \leq \frac12
\epsilon \|\phi^0(t)\|^2_1 \leq \frac12\epsilon B^2(t)
\|\phi_0\|^2_1.
\end{equation}
\end{lemma}

\noindent \it Remark. \rm Note that $\phi^0(t) = e^{iLt}\phi_0$ by
definition. Assumption \eqref{Lcon2} says that this unitary
evolution is bounded in the $H^1(\Gamma)$ norm.

\begin{proof}
The regularity guaranteed by conditions \eqref{Lcon1},
\eqref{Lcon2} and Theorem~\ref{exun} allows us to multiply the
equation
\[
(\phi^\epsilon - \phi^0)' = iL(\phi^\epsilon -
\phi^0) - \epsilon \Gamma \phi^\epsilon
\]
by $\phi^\epsilon - \phi^0.$
We obtain
\[
\frac d{dt} \|\phi^\epsilon - \phi^0\|^2  \leq 2\epsilon
(\|\phi^\epsilon\|_1 \|\phi^0\|_1 - \|\phi^\epsilon\|^2_1) \le
\frac12\epsilon \|\phi^0 \|_1^2,
\]
which is the first inequality in \eqref{l2approx}. The second
inequality follows simply from the assumption~\eqref{Lcon2}.
\end{proof}

The following corollary is immediate.
\begin{corollary}\label{expapp}
Assume that \eqref{Lcon1} and \eqref{Lcon2} are satisfied, and the
initial data $\phi_0 \in H^1.$ Then the solutions
$\phi^\epsilon(t)$ and $\phi^0(t)$ defined in Lemma~\ref{approxH}
satisfy
\[ \|\phi^\epsilon(t) - \phi^0(t)\|^2 \leq  \frac12\epsilon
\|\phi_0\|_1^2 \int_0^\tau B^2(t)\,dt\]
for any time $t \leq \tau.$
\end{corollary}

Finally, we observe that conditions \eqref{Lcon1} and
\eqref{Lcon2} are independent. Taking $L=\Gamma$ shows that
\eqref{Lcon2} does not imply \eqref{Lcon1}, because in this case
the evolution $e^{iLt}$ is unitary on $H^1$ but the domain of $L$ is
$H^2\subsetneq H^1$. On the other hand, \eqref{Lcon1} does not imply
\eqref{Lcon2}, as is the case in the following example. Let $H\equiv
L^2(0,1)$, define the operator $\Gamma$ by $\Gamma f(x) \equiv
\sum_n e^{n^2}\hat f(n)e^{2\pi inx}$ for all $f\in H$ such that
$e^{n^2}\hat f(n)\in \ell^2(\Zm)$, and take $Lf(x)\equiv xf(x)$.
Then $L$ is bounded on $H$ and so \eqref{Lcon1} holds automatically, but
\[
\left[e^{itL}f\right](x)=f(x)e^{itx}
\]
so that
$e^{
2\pi iL} e^{2\pi inx}=e^{2\pi i(n+1)x}$. It follows that $e^{2\pi iL}$ is
not bounded on $H^1$ (and neither is $e^{iLt}$ for any $t\neq 0$).

\section{The Abstract Criterion}\label{abs}

One direction in the proof of Theorem~\ref{general} is much easier. We
start by proving this easy direction: that existence of $H^1(\Gamma)$
eigenvectors of $L$ ensures existence of $\tau,\delta>0$ and $\phi_0$
with $\|\phi_0\|=1$ such that $\|\phi^A(\tau)\| > \delta$ for all $A$
-- that is, if such eigenvectors exist, then the
operator $L$ is not relaxation enhancing.

\begin{proof}[Proof of the first part of Theorem~\ref{general}]
Assume that the initial data $\phi_0 \in H^1$ for (\ref{hseq}) is an
eigenvector of $L$ corresponding to an eigenvalue $E$, normalized so that
$\|\phi_0\|=1$. Take the inner
product of (\ref{hseq}) with $\phi_0.$ We arrive at
\[
\frac{d}{dt} \langle \phi^A(t), \phi_0 \rangle = iA E \langle
\phi^A(t), \phi_0 \rangle - \langle \Gamma \phi^A(t), \phi_0
\rangle.
\]
This and the assumption $\phi_0 \in H^1$ lead to
\[
\left| \frac{d}{dt} \left( e^{-i A E t} \langle \phi^A(t), \phi_0
\rangle \right) \right| \leq \frac12 \left( \| \phi^A(t) \|_1^2 +
\|\phi_0\|^2_1 \right).
\]
Note that the value of the expression being differentiated on the
left hand side is equal to one at $t=0.$ By Lemma~\ref{gradbound}
(with a simple time rescaling) we have $\int_0^\infty \|\phi^A(t)\|_1^2
\,dt \leq 1/2.$ Therefore, for $t \leq \tau =
(2\|\phi_0\|_1^2)^{-1}$ we have $|\langle \phi^A(t), \phi_0 \rangle|
\geq 1/2.$ Thus, $\|\phi^A(\tau)\| \geq 1/2,$ uniformly in $A.$
\end{proof}
Note also that we have proved that in the presence of an $H^1$
eigenvector of $L$, enhanced relaxation does not happen for some $\phi_0$
even if we allow $A(\tau,\delta)$ to be $\phi_0$-dependent as well.
This explains Remark 2 after Definition \ref{relen}.

The proof of the converse is more subtle, and will require some
preparation. We switch to the equivalent formulation \eqref{sdeq}.  We
need to show that if $L$ has no $H^1$ eigenvectors, then for all
$\tau, \delta>0$ there exists $\epsilon_0(\tau,\delta)>0$ such that if
$\epsilon < \epsilon_0,$ then $\|\phi^\epsilon(\tau/\epsilon)\| <
\delta$ whenever $\|\phi_0\|=1$. The main idea of the proof can be
naively described as follows.  If the operator $L$ has purely
continuous spectrum or its eigenfunctions are rough then the
$H^1$-norm of the free evolution (with $\epsilon=0$) is large most of
the time. However, the mechanism of this effect is quite different for
the continuous and point spectra. On the other hand, we will show
that for small $\epsilon$ the full evolution is close to the free
evolution for a sufficiently long time. This clearly leads to
dissipation enhancement.

The first ingredient that we need to recall is the so-called RAGE
theorem.
\begin{thm}[RAGE] \label{RAGEthm}
Let $L$ be a self-adjoint operator in a Hilbert space $H.$ Let
$P_c$ be the spectral projection on its continuous spectral
subspace. Let $C$ be any compact operator. Then for any $\phi_0
\in H,$ we have
\[ \lim_{T \rightarrow \infty} \frac1T \int\limits_0^T \|C
e^{iLt} P_c \phi_0\|^2\,dt = 0. \]
\end{thm}

Clearly, the result can be equivalently stated for a unitary operator
$U$, replacing $e^{iLt}$ with $U^t.$ The proof of the RAGE theorem can
be found, for example, in \cite{CFKS}.

A direct consequence of the RAGE theorem is the following lemma.
Recall that we denote by $0<\lambda_1\le\lambda_2\le\dots$ the
eigenvalues of the operator $\Gamma$ and by $e_1,e_2,\dots$ the
corresponding orthonormal eigenvectors. Let us also denote by $P_N$
the orthogonal projection on the subspace spanned by the first $N$ eigenvectors
$e_1, \dots, e_N$ and by $S= \{ \phi\in H:~
\|\phi\|=1\}$ the unit sphere in $H.$ The following lemma shows that
if the initial data lies in the continuous spectrum of $L$ then the $L$-evolution
will spend most of time in the higher modes of $\Gamma$.
\begin{lemma}\label{RAGE}
Let $K \subset S$ be a compact set. For any $N,\sigma>0,$ there
exists $T_c(N,\sigma,K)$ such that for all $T \geq T_c(N,\sigma, K)$
and any $\phi \in K,$ we have
\begin{equation}\label{RAGEbound}
\frac1T \int\limits_0^T \|P_N e^{iLt}P_c\phi\|^2\,dt \leq \sigma.
\end{equation}
\end{lemma}

\noindent \it Remark. \rm The key observation of Lemma \ref{RAGE}
is that the time $T_c(N,\sigma,K)$ is uniform for all $\phi \in
K.$

\begin{proof}
Since $P_N$ is compact, we see that for any vector $\phi \in S,$ there
exists a time $T_c(N,\sigma, \phi)$ that depends on the function
$\phi$ such that \eqref{RAGEbound} holds for $T>T_c(N,\sigma,\phi)$ --
this is assured by Theorem \ref{RAGEthm}. To prove the uniformity in
$\phi,$ note that the function
\[
f(T, \phi) = \frac1T \int\limits_0^T \|P_N e^{iLt}P_c\phi\|^2\,dt
\]
is uniformly continuous on $S$ for all $T$ (with constants
independent of $T$):
\begin{eqnarray*}
&&\!\!\!\!|f(T,\phi)-f(T,\psi)|\le\frac1T \int\limits_0^T \left|
\|P_N e^{iLt}P_c\phi\|-\|P_N e^{iLt}P_c\psi\|\right|\left(
\|P_N e^{iLt}P_c\phi\|+\|P_N e^{iLt}P_c\psi\|\right)\,dt\\
&&\,\,\,\,\,\,\,\,\,\,\,\,\,\,\,\,\,\,\,\,\,\,\,\,
\,\,\,\,\,\,\,\,\,\,\,\,\,\,\,\,\,\,\,\,\,\le
(\|\phi\|+\|\psi\|)\frac1T \int\limits_0^T \|P_N
e^{iLt}P_c(\phi-\psi)\|dt\le 2\|\phi-\psi\|.
\end{eqnarray*}
Now, existence of a uniform $T_c(N,\sigma, K)$
follows from compactness of $K$ by standard arguments.
\end{proof}

We also need a lemma which controls from below the growth of the
$H^1$ norm of free solutions corresponding to rough
eigenfunctions. We denote by $P_p$ the spectral projection on the
pure point spectrum of the operator $L$.
\begin{lemma}\label{roughgrowth}
Assume that not a single eigenvector of the operator $L$ belongs to
$H^1(\Gamma).$ Let $K \subset S$ be a compact set. Consider the set
$K_1\equiv \{\phi\in K \,|\, \|P_p\phi\| \geq 1/2\}.$ Then for any
$B>0$ we can find $N_p(B,K)$ and $T_p(B,K)$ such that for any $N\ge
N_p(B,K)$, any $T \geq T_p(B,K)$ and any $\phi \in K_1$, we have
\begin{equation}\label{h1normgr}
\frac{1}{T} \int\limits_0^T \| P_{N} e^{iLt}P_p\phi\|_1^2\,dt \geq
B.
\end{equation}
\end{lemma}
\noindent {\it Remark.} Note that unlike in \eqref{RAGEbound}, we
have the $H^1$ norm in \eqref{h1normgr}.

\begin{proof}
The set $K_1$ may be empty, in which case there is nothing to
prove. Otherwise, let us denote by $E_j$ the eigenvalues of $L$ (distinct,
without repetitions) and by $Q_j$ the orthogonal projection on the
space spanned by the eigenfunctions corresponding to $E_j.$ First, let
us show that for any $B>0$ there is $N(B,K)$ such that for any
$\phi \in K_1$
we have
\begin{equation}\label{bigh1}
\sum\limits_{j} \|P_N Q_j \phi\|^2_1 \geq 2B
\end{equation}
if $N \geq N(B,K).$ It is clear that for each fixed $\phi$ with $P_p
\phi \ne 0$ we can find $N(B,\phi)$ so that \eqref{bigh1} holds,
since by assumption $Q_j \phi$ does not belong to $H^1$ whenever
$Q_j\phi \ne 0$. Assume that $N(B,K)$ cannot be chosen uniformly for
$\phi \in K_1.$ This means that for any $n,$ there exists $\phi_n
\in K_1$ such that
\[
\sum\limits_{j} \|P_n Q_j \phi_n\|^2_1 < 2B.
\]
Since $K_1$ is compact, we can find a subsequence $n_l$ such that
$\phi_{n_l}$ converges to $\tilde{\phi}\in K_1$ in $H$ as $n_l
\rightarrow \infty.$ For any $N$ and any $n_{l_1}
>N$ we have
\[
\sum_{j} \|P_N Q_j \tilde{\phi}\|^2_1 \leq \sum_{j} \|P_{n_{l_1}}
Q_j \tilde{\phi}\|^2_1 \le \liminf_{l \rightarrow \infty}\sum_{j}
\|P_{n_{l_1}} Q_j \phi_{n_l}\|^2_1.
\]
The last inequality follows by Fatou's Lemma from the convergence
of $\phi_{n_l}$ to $\tilde{\phi}$ in $H$ and the fact that
$\|P_{n_{l_1}} Q_j \psi \|_1 \leq \lambda_{n_{l_1}}^{1/2}\|\psi\|$
for any $n_{l_1}.$ But now the expression on the right hand side
is less than or equal to
\[
\liminf_{l\to\infty} \sum_{j} \|P_{n_{l}} Q_j \phi_{n_l}\|^2_1 \le
2B.
\]
Thus $\sum\limits_{j}
 \|P_N Q_j \tilde{\phi}\|^2_1 \leq 2B$ for any $N,$ a
contradiction since $\tilde{\phi} \in K_1.$

Next, take $\phi \in K_1$ and consider
\begin{equation}\label{dyn11}
\frac1T \int\limits_0^T \|P_N e^{iLt} P_p\phi \|_1^2\,dt =
\sum\limits_{j,l} \frac{e^{i(E_j - E_l)T} -1}{i(E_j - E_l)T}
\langle \Gamma P_N Q_j \phi, P_N Q_l \phi \rangle.
\end{equation}
In \eqref{dyn11}, we set $(e^{i(E_j - E_l)T} -1)/i(E_j - E_l)T\equiv
1$ if $j=l.$ Notice that the
sum above converges
absolutely. Indeed, $P_N Q_j \phi = \sum_{i=1}^N
\langle Q_j \phi, e_i \rangle e_i,$ and $\langle \Gamma e_i, e_k
\rangle = \lambda_i \delta_{ik},$ therefore
\[
\langle \Gamma P_N Q_j \phi, P_N Q_l \phi \rangle =
\sum\limits_{i=1}^N \lambda_i \langle Q_j \phi, e_i \rangle
\overline{\langle Q_l\phi, e_i \rangle} .
\]
Therefore the sum on the right hand side of \eqref{dyn11} does not
exceed
\begin{align}
\sum\limits_{i=1}^N \lambda_i \sum\limits_{j,l} |\langle Q_j \phi,
e_i \rangle \langle Q_l \phi, e_i \rangle | & \leq \lambda_N
\sum\limits_{i=1}^N \sum\limits_{j,l}
\|Q_j\phi\|\|Q_l\phi\||\langle Q_j \phi/\|Q_j\phi\|, e_i \rangle
\langle Q_l\phi/\|Q_l\phi\| , e_i \rangle | \notag
\\ & \le \lambda_N \sum\limits_{i=1}^N \sum_{j,l} \|Q_l \phi\|^2 |\langle
Q_j\phi/\|Q_j\phi\|, e_i \rangle|^2  \leq  \lambda_N N ,
\label{absolute}
\end{align}
with the second step obtained using the Cauchy-Schwartz inequality,
and the third by
$\|\phi\|=\|e_i\|=1$. Then for each fixed $N,$ we have by the
dominated convergence theorem that the expression in \eqref{dyn11}
converges to $\sum_j \|\Gamma^{1/2} P_N Q_j \phi\|^2=\sum_j \| P_N
Q_j \phi\|_1^2$ as $T \rightarrow \infty.$ Now assume $N \geq
N_p(B,K)\equiv N(B,K),$ so that \eqref{bigh1} holds. We claim that
we can choose $T_p(B,K)$ so that for any $T\ge T_p(B,K)$ we have
\begin{equation}\label{offdiag}
\left| \frac 1T \int_0^T \| P_N e^{iLt} P_p\phi \|_1^2 - \sum_j \|
P_N Q_j \phi\|_1^2 \right| = \left| \sum\limits_{l \ne j}
\frac{e^{i(E_j - E_l)T} -1}{i(E_j - E_l)T} \langle \Gamma P_N Q_j
\phi, P_N Q_l \phi \rangle \right| \leq B
\end{equation}
for all $\phi \in K_1.$ Indeed, this follows from convergence to
zero for each individual $\phi$ as $T \rightarrow \infty,$
compactness of $K_1,$ and uniform continuity of the expression in
the middle of \eqref{offdiag} in $\phi$ for each $T$ (with constants
independent of $T$). The latter is proved by estimating the
difference of these expressions for $\phi,\psi\in K_1$ and any $T$
by
\[
\sum\limits_{l \ne j} |\langle \Gamma P_N Q_j \phi, P_N Q_l
(\phi-\psi) \rangle | + |\langle \Gamma P_N Q_j (\phi-\psi), P_N
Q_l \psi \rangle |,
\]
which is then bounded by $2\lambda_N N \|\phi-\psi\|$ using the
trick from \eqref{absolute}. Combining \eqref{bigh1} and
\eqref{offdiag} proves the lemma.
\end{proof}

We can now complete the proof of Theorem~\ref{general}.

\begin{proof}[Proof of Theorem~\ref{general}]
Recall that given any $\tau,\delta >0,$ we need to show the
existence of $\epsilon_0>0$ such that if $\epsilon < \epsilon_0,$
then $\| \phi^\epsilon (\tau/\epsilon)\| < \delta$ for any initial
datum $\phi_0 \in H,$ $\|\phi_0\|= 1.$ Here $\phi^\epsilon(t)$ is
the solution of \eqref{sdeq}. Let us outline the idea of the
proof. Lemma~\ref{easy} tells us that if the $H^1$ norm of the
solution $\phi^\epsilon(t)$ is large, relaxation is happening
quickly. If, on the other hand, $\|\phi^\epsilon(\tau_0)\|_1^2
\leq \lambda_M \|\phi^\epsilon(\tau_0)\|^2,$ where $M$ is to be
chosen depending on $\tau, \delta$, then the set of all unit
vectors satisfying this inequality is compact, and so we can apply
Lemma~\ref{RAGE} and Lemma~\ref{roughgrowth}. Using these lemmas,
we will show that even if the $H^1$ norm is small at some moment
of time $\tau_0$, it will be large on the average in some time
interval after $\tau_0.$ Enhanced relaxation will follow.

We now provide the details. Since $\Gamma$ is an unbounded
positive operator with a discrete spectrum, we know that its
eigenvalues $\lambda_n \rightarrow \infty$ as $n \rightarrow
\infty.$  Let us choose $M$ large enough, so that $e^{-\lambda_M
\tau/80} < \delta.$ Define the sets $K\equiv\{ \phi \in S \,|\,
\|\phi\|_1^2 \leq \lambda_M \}\subset S$ and as before, $K_1\equiv
\{\phi\in K\,|\, \|P_p\phi\| \geq 1/2\}$. It is easy to see that
$K$ is compact. Choose $N$ so that $N \geq M$ and $N\ge
N_p(5\lambda_M,K)$ from Lemma \ref{roughgrowth}. Define
\[
\tau_1 \equiv {\rm max} \left\{ T_p(5\lambda_M,K), T_c(N, \tfrac
 {\lambda_M}{20 \lambda_N},K) \right\},
\]
where $T_p$ is from Lemma~\ref{roughgrowth}, and $T_c$ from
Lemma~\ref{RAGE}. Finally, choose $\epsilon_0>0$ so that $\tau_1<
\tau/2 \epsilon_0,$ and
\begin{equation}\label{eps0cond4}
\epsilon_0 \int\limits_0^{\tau_1} B^2(t)\,dt \leq \frac{1}{20
\lambda_N},
\end{equation}
where $B(t)$ is the function from condition \eqref{Lcon2}.

Take any $\epsilon < \epsilon_0.$ If we have  $\|\phi^\epsilon (s)\|_1^2 \geq
\lambda_M \|\phi^\epsilon(s)\|^2$ for all $s\in[0,\tau]$ then Lemma  \ref{easy}
implies that $\|\phi^\epsilon(\tau/\epsilon)\|\le
e^{-2\lambda_M\tau}\le\delta$ by the choice of $M$ and we are
done. Otherwise, let $\tau_0$ be the first time in the interval
$[0,\tau/\epsilon]$ such that $\|\phi^\epsilon (\tau_0)\|_1^2 \leq
\lambda_M \|\phi^\epsilon(\tau_0)\|^2$ (it may be that $\tau_0=0,$
of course). We claim that the following estimate holds for the decay
of $\|\phi^\epsilon(t)\|$
on the interval $[\tau_0,\tau_0+\tau_1]$:
\begin{equation}\label{ccase5-1}
\|\phi^\epsilon(\tau_0+\tau_1)\|^2
\leq e^{-\lambda_M \epsilon \tau_1/20}\|\phi^\epsilon(\tau_0)\|^2.
\end{equation}

For the sake of transparency, henceforth we will denote
$\phi^\epsilon(\tau_0)=\phi_0.$ On the interval
$[\tau_0,\tau_0+\tau_1],$ consider the function $\phi^0(t)$
satisfying $\tfrac d{dt}\phi^0(t) = iL\phi^0(t),$ $\phi^0(\tau_0) =
\phi_0.$ Note that by the choice of $\epsilon_0$, \eqref{eps0cond4}
and Corollary~\ref{expapp}, we have
\begin{equation}\label{close11}
\| \phi^\epsilon(t) - \phi^0(t) \|^2 \leq \frac{\lambda_M}{40
\lambda_N} \|\phi_0\|^2
\end{equation}
for all $t \in [\tau_0, \tau_0+\tau_1].$ Split $\phi^0(t)=
\phi_c(t)+\phi_p(t),$ where $\phi_{c,p}$ also solve the free
equation $\tfrac d{dt}\phi_{c,p}(t) = iL\phi_{c,p}(t),$ but with
initial data $P_c\phi_0$ and $P_p \phi_0$ at $t=\tau_0,$
respectively. We will now consider two cases.

{\it Case I.} Assume that $\|P_c \phi_0\|^2 \geq \frac34
\|\phi_0\|^2,$ or, equivalently, $\|P_p \phi_0\|^2 \leq \frac14
\|\phi_0\|^2.$ Note that since $\phi_0 /\|\phi_0\| \in K$ by the
hypothesis, we can apply Lemma~\ref{RAGE}. Our choice of $\tau_1$
implies that
\begin{equation}\label{ccase1}
\frac{1}{\tau_1} \int\limits_{\tau_0}^{\tau_0+\tau_1} \|P_N
\phi_c(t)\|^2 \,dt \leq \frac{\lambda_M}{20 \lambda_N}
\|\phi_0\|^2.
\end{equation}
By elementary considerations,
\[
\!\|(I-P_N)\phi^0(t)\|^2 \geq \frac 12 \|(I-P_N) \phi_c(t)\|^2 - \|
(I-P_N)\phi_p(t)\|^2 \ge \frac 12\|\phi_c(t)\|^2 - \frac
12\|P_N\phi_c(t)\|^2 - \|\phi_p(t)\|^2.
\]
Taking into account the fact that the free evolution $e^{iLt}$ is
unitary, $\lambda_N\ge\lambda_M$, our assumptions on
$\|P_{c,p}\phi_0\|$ and \eqref{ccase1}, we obtain
\begin{equation}\label{ccase2}
\frac{1}{\tau_1} \int\limits_{\tau_0}^{\tau_0+\tau_1} \|(I-P_N)
\phi^0(t)\|^2 \,dt \geq \frac{1}{10} \|\phi_0\|^2.
\end{equation}
Using \eqref{close11}, we conclude that
\begin{equation}\label{ccase3}
\frac 1{\tau_1}\int\limits_{\tau_0}^{\tau_0+\tau_1} \|(I-P_N)
\phi^\epsilon(t)\|^2 \,dt \geq \frac{1}{40} \|\phi_0\|^2.
\end{equation}
This estimate implies that
\begin{equation}\label{ccase4}
\int\limits_{\tau_0}^{\tau_0+\tau_1} \|\phi^\epsilon(t)\|^2_1 \,dt
\geq \frac{\lambda_N \tau_1}{40} \|\phi_0\|^2.
\end{equation}
Combining \eqref{ccase4} with \eqref{dissip} yields
\begin{equation}\label{ccase5}
\|\phi^\epsilon(\tau_0+\tau_1)\|^2 \le \left( 1-\frac
{\lambda_N\epsilon\tau_1}{20} \right)\|\phi^\epsilon(\tau_0)\|^2
\leq e^{-\lambda_N \epsilon \tau_1/20}\|\phi^\epsilon(\tau_0)\|^2.
\end{equation}
This finishes the proof of (\ref{ccase5-1}) in the first case since $\lambda_N\ge\lambda_M$.

{\it Case II.} Now suppose that $\|P_p \phi_0\|^2 \geq \frac14
\|\phi_0\|^2.$ In this case $\phi_0/\|\phi_0\| \in K_1,$ and we
can apply Lemma~\ref{roughgrowth}. In particular, by the choice of
$N$ and $\tau_1,$ we have
\begin{equation}\label{pcase1}
\frac{1}{\tau_1} \int\limits_{\tau_0}^{\tau_0+\tau_1} \|P_N
\phi_p(t)\|^2_1 \,dt \geq 5 \lambda_M \|\phi_0\|^2.
\end{equation}
Since \eqref{ccase1} still holds because of our choice of $\tau_0$ and $\tau_1$, it follows
that
\begin{equation}\label{pcase2}
\frac{1}{\tau_1} \int\limits_{\tau_0}^{\tau_0+\tau_1} \|P_N
\phi_c(t)\|^2_1 \,dt \leq  \frac{\lambda_M}{20} \|\phi_0\|^2.
\end{equation}
Note that the $H$-norm in (\ref{ccase1}) has been replaced in
(\ref{pcase2}) by the $H^1$-norm at the expense of the factor of
$\lambda_N$.  Together, \eqref{pcase1} and \eqref{pcase2} imply
\begin{equation}\label{pcase3}
\frac{1}{\tau_1} \int\limits_{\tau_0}^{\tau_0+\tau_1} \|P_N
\phi^0(t)\|^2_1 \,dt \geq 2 \lambda_M \|\phi_0\|^2.
\end{equation}
Finally, \eqref{pcase3} and \eqref{close11} give
\begin{equation}\label{pcase4}
\int\limits_{\tau_0}^{\tau_0+\tau_1} \|P_N \phi^\epsilon(t)\|^2_1
\,dt \geq  \frac{\lambda_M \tau_1}2 \|\phi_0\|^2
\end{equation}
since $\|P_N\phi^\epsilon-P_N\phi^0\|_1^2\le \lambda_N\|\phi^\epsilon-\phi^0\|^2$.
As before, \eqref{pcase4} implies
\begin{equation}\label{pcase5}
\|\phi^\epsilon(\tau_0+\tau_1)\|^2 \leq e^{-\lambda_M \epsilon
\tau_1}\|\phi^\epsilon(\tau_0)\|^2,
\end{equation}
finishing the proof of (\ref{ccase5-1}) in the second case.

Summarizing, we
see that if $\|\phi^\epsilon (\tau_0)\|_1^2 \leq \lambda_M
\|\phi^\epsilon(\tau_0)\|^2,$ then
\begin{equation}\label{final1}
\|\phi^\epsilon(\tau_0+\tau_1)\|^2 \leq e^{-\lambda_M \epsilon
\tau_1/20}\|\phi^\epsilon(\tau_0)\|^2.
\end{equation}
On the other hand, for any interval $I=[a,b]$ such that
$\|\phi^\epsilon(t)\|^2_1 \geq \lambda_M \|\phi^\epsilon(t)\|^2$
on $I,$ we have by Lemma~\ref{easy} that
\begin{equation}\label{final2}
\|\phi^\epsilon(b)\|^2 \leq e^{-2\lambda_M \epsilon
(b-a)}\|\phi^\epsilon(a)\|^2.
\end{equation}
Combining all the decay factors gained from \eqref{final1} and
\eqref{final2}, and using $\tau_1<\tau/2\epsilon$, we find that
there is $\tau_2\in[\tau/2\epsilon,\tau/\epsilon]$ such that
\[
\|\phi^\epsilon(\tau_2)\|^2 \leq e^{-\lambda_M
\epsilon\tau_2/20}\leq e^{-\lambda_M \tau/40} < \delta^2
\]
by our choice of $M.$ Then \eqref{dissip} gives
$\|\phi^\epsilon(\tau/\epsilon)\|\le
\|\phi^\epsilon(\tau_2)\|<\delta$, finishing the proof of
Theorem~\ref{general}.
\end{proof}

\section{Examples With Rough Eigenvectors}\label{rough}

It is not immediately obvious that condition \eqref{Lcon2},
$\|e^{iLt}\phi\|_1 \leq B(t) \|\phi\|_1$ for any $\phi_0 \in H^1,$
is consistent with the existence of eigenvectors of $L$ which are
not in $H^1.$ The purpose of this section is to show that, in
general, rough eigenvectors may indeed be present under the
conditions of Theorem~\ref{general}. We provide here two simple
examples of operators $\Gamma$ and $L$ in which \eqref{Lcon2} is
satisfied and $L$ has only rough eigenfunctions. In both cases $L$
will be a discrete Schr\" odinger operator on $\Zm^+$ resp., more
generally, a Jacobi matrix, and $\Gamma$ a multiplication
operator. One more example with rough eigenfunctions will deal
with an actual fluid flow and will be discussed in
Section~\ref{flowex}.

The first is an explicit example with one rough eigenvector that
is a discrete version of the celebrated Wigner-von Neumann
construction \cite{WvN} of an imbedded eigenvalue of a
Schr\"odinger operator with a decaying potential. The second
example is implicit, its existence being guaranteed by a result of
Killip and Simon \cite{KS}, and demonstrates that all eigenvectors
of $L$ can be rough while at the same time the eigenvalues can be
dense in the spectrum of $L$.
\smallskip

{\it Example 1.} Let $\tilde{\Gamma}$ be the operator of
multiplication by $n$ on $l^2(\Zm^+),$ $\Zm^+=\{1,2,\dots,\}.$
Furthermore, let $\tilde{L}$ be the discrete Schr\"odinger operator
on $l^2(\Zm^+)$
\[
\tilde{L}u_n = u_{n+1}+u_{n-1}+v_n u_n
\]
for $n\ge 1$, with the potential
\[
v_n\equiv \begin{cases} -\tfrac 2{n+2} & n \text{ even},
\\ \tfrac 2{n-1} & n>1 \text{ odd},
\\ -1 & n=1,
\end{cases}
\]
and the self-adjoint boundary condition $u_0\equiv 0$. Then
$\tilde L$ has eigenvalue zero with eigenfunction $u$ given by
\[
u_{2n-1}=u_{2n}= \tfrac{(-1)^n}n
\]
for $n\ge 1$, because then $\tilde Lu\equiv 0$ and
$u\in\ell^2(\Zm^+)$. Note that $u$ does not belong to
$H^1(\tilde{\Gamma}).$

It is not difficult to show that $\tilde L$ has no more
eigenvalues in its essential spectrum $[-2,2]$ (for example, using
the so-called EFGP transform, see \cite{KLS} for more details).
The eigenvalue zero is a consequence of a resonant structure of
the potential which is tuned to this energy. There may be (and
there are) other eigenvalues outside $[-2,2],$ with eigenfunctions
that are exponentially decaying and so do belong to
$H^1(\tilde{\Gamma}).$ It is also known that $\tilde L$ has no
singular continuous spectrum and it has absolutely continuous
spectrum that fills $[-2,2]$. More precisely, the absolutely
continuous part of the spectral measure gives positive weight to
any set of positive Lebesgue measure lying in $[-2,2]$ (see, e.g.,
\cite{KS}).

To get an example where we have only rough eigenfunctions, we will
project away the eigenfunctions lying in $H^1.$ Namely, denote by
$D$ the subspace spanned by all eigenfunctions of $\tilde{L}$,
with the exception of $u.$ Denote $P$ the projection on the
orthogonal complement of $D,$ and set $\Gamma = P \tilde{\Gamma}
P,$ $L = P\tilde{L}P.$ Then $\Gamma,$ $L$ are self-adjoint on the
infinite dimensional Hilbert space $H=Pl^2(\Zm^+)$, and by
construction $L$ has absolutely continuous spectrum filling
$[-2,2]$ as well as a single eigenvalue equal to zero. The
corresponding eigenfunction is $u$ and it does not belong to
$H^1(\Gamma)$ because
\[
\langle \Gamma u, u \rangle = \langle P \tilde{\Gamma} P u,u \rangle
= \langle \tilde{\Gamma}u,u \rangle \ge \sum_n |n|(n^{-1})^2 =
\infty.
\]

Let us check the conditions of Theorem \ref{general}. First,
$\Gamma$ is positive because $\tilde\Gamma$ is. It is also unbounded
and has a discrete spectrum. Indeed, let $H_R \subset H$ be the
subspace of all vectors $\phi \in H$ such that
\begin{equation}\label{dse}
\langle \Gamma \phi, \phi \rangle \leq R\langle \phi, \phi \rangle.
\end{equation}
Then for each such $\phi$ we also have \eqref{dse} with
$\tilde{\Gamma}$ instead of $\Gamma.$ By the minimax principle for
self-adjoint operators this implies that $\lambda_n \geq
\tilde{\lambda}_n=n,$ where $\lambda_n,$ $\tilde{\lambda}_n$ are
the $n$-th eigenvalues of $\Gamma$ and $\tilde{\Gamma},$
respectively (counting multiplicities).

Also, $L$ is a bounded operator on $H$ (since $\tilde{L}$ is) and so
\eqref{Lcon1} is satisfied automatically. Finally, observe that for
any $\phi \in H^1(\Gamma),$ we have
\begin{equation}\label{dse1}
|\langle \Gamma L \phi, \phi \rangle| = |\langle P \tilde{\Gamma} P
\tilde{L} P \phi,  \phi \rangle| = | \langle \tilde{\Gamma}
\tilde{L} \phi, \phi \rangle | \leq C\|\phi\|_1^2
\end{equation}
The second equality in \eqref{dse1} follows from the fact that
$\tilde{L}$ and $P$ commute by construction and $P\phi=\phi$ for
$\phi \in H.$ The inequality in \eqref{dse1} holds since $\|\tilde
L\phi\|_1\le C\|\phi\|_1$, which follows from the fact that
$\tilde{L}$ is tridiagonal and both
$\tilde\lambda_{n+1}/\tilde\lambda_n$ and $v_n$ are bounded. Now
given $\phi \in H^1(\Gamma),$ set $\phi(t) = e^{iLt}\phi.$ Then
\[
\left| \frac{d}{dt} \|\phi(t)\|_1^2 \right| \leq 2|\langle \Gamma L
\phi(t), \phi(t) \rangle| \leq C \|\phi(t)\|_1^2
\]
by \eqref{dse1}. This a priori estimate and Gronwall's Lemma allow
one to conclude that \eqref{Lcon2} holds with $B(t) = e^{Ct/2}.$
This concludes our first example.
\smallskip

{\it Example 2.} We let $H\equiv \ell^2(\Zm^+)$ and define $\Gamma$
to be the multiplication by $e^n$. In order to provide an example
with a much richer set of rough eigenfunctions, we will now consider
$L$ to be a Jacobi matrix
\[
Lu_n = a_{n}u_{n+1}+a_{n-1}u_{n-1}+v_n u_n,
\]
with $a_n>0$, $v_n\in\mathbb{R}$ and boundary condition $u_0\equiv
0$. We choose $\nu$ to be a pure point measure of total mass
$\tfrac 12$, whose mass points are contained and dense in
$(-2,2)$, and define the probability measure $d\mu(x)\equiv
d\nu(x) + \tfrac 18\chi_{[-2,2]}(x)dx$. By the Killip-Simon
\cite{KS} characterization of spectral measures of Jacobi matrices
that are Hilbert-Schmidt perturbations of the free half-line
Schr\"odinger operator (with $a_n=1$, $v_n=0$), there is a unique
Jacobi matrix $L$ such that $a_n-1,v_n\in\ell^2(\Zm^+)$ and its
spectral measure is $\mu$. In particular, the eigenvalues of $L$
are dense in its spectrum $\sigma(L)=[-2,2]$.

The conditions of Theorem \ref{general} are again satisfied, with
the key estimate $\|L\phi\|_1\le C\|\phi\|_1$ holding because
$\lambda_{n+1}/\lambda_n,a_n,v_n$ are bounded. Moreover, it is
easy to show (see below) that the fact that eigenvalues of $L$ are
inside $(-2,2)$ and $a_n-1,v_n\in\ell^2$ imply that eigenfunctions
of $L$ decay slower than $e^{-C\sqrt n}$ for some $C$. More
precisely, if $u$ is an eigenfunction of $L$, then $\lim_n
(u_n^2+u_{n-1}^2)e^{C\sqrt n}=\infty$, and so obviously $u\notin
H^1(\Gamma)$ (actually, $u\notin H^s(\Gamma)$ for any $s>0$).

To obtain the well-known bound on the eigenfunction decay, let $u$
be an eigenfunction of $L$ corresponding to eigenvalue $E\in
(-2,2)$, that is,
\begin{equation}\label{efunc}
Eu_n=a_{n}u_{n+1}+a_{n-1}u_{n-1}+v_n u_n
\end{equation}
for $n\ge 1$. Define the square of the Pr\" ufer amplitude of $u$
by
\[
R_n\equiv u_n^2 + u_{n-1}^2 - Eu_nu_{n-1} = \frac{2-|E|}2 (u_n^2 +
u_{n-1}^2) + \frac {|E|}2 (u_n-u_{n-1})^2>0
\]
and $c_n\equiv |a_n-1|+|a_{n-1}-1|+|v_n|\in\ell^2$. After
expressing $u_{n+1}$ in terms of $u_n$ and $u_{n-1}$ using
\eqref{efunc}, one obtains (with each $|O(c_n)|\le C_Ec_n$)
\[
\frac{R_{n+1}}{R_n} = (1+O(c_n)) \frac {\frac{2-|E|}2
((1+O(c_n))u_n^2 + (1+O(c_n))u_{n-1}^2) + \frac {|E|}2 (1+O(c_n))
(u_n-u_{n-1})^2} {\frac{2-|E|}2 (u_n^2 + u_{n-1}^2) + \frac {|E|}2
(u_n-u_{n-1})^2}
\]
if $E\neq 0$ and
\[
\frac{R_{n+1}}{R_n} = (1+O(c_n)) \frac {(1+O(c_n))u_n^2 +
(1+O(c_n))u_{n-1}^2 + O(c_n)u_nu_{n-1}} {u_n^2 + u_{n-1}^2}
\]
if $E=0$. In either case, $R_{n+1}/R_n=1+O(c_n)$, which means that
\[
R_n\ge R_{n_0}\prod_{k=n_0+1}^n(1-C_Ec_k)\ge
R_{n_0}\exp(-2C_E\sum_{k=n_0+1}^n c_k)\ge
R_{n_0}\exp(-2C_E\|c_k\|_2\sqrt{n})
\]
if $n_0$ is chosen so that $C_Ec_k<\tfrac 12$ for $k>n_0$. But
then the definition of $R_n$ shows that $\lim_n
(u_n^2+u_{n-1}^2)e^{C\sqrt n}=\infty$ for some $C<\infty$. This
concludes the example.
\smallskip

We have thus proved

\begin{thm}
There exist a self-adjoint, positive, unbounded operator $\Gamma$
with a discrete spectrum and a self-adjoint operator $L$ such that
the following conditions are satisfied.
\begin{itemize}
\item $\|L\phi\|\leq C\|\phi\|_1$ and $\|e^{iLt}\phi\|_1\leq
B(t)\|\phi\|_1$ for some $C<\infty$, $B(t)\in L^1_{\rm
loc}[0,\infty)$ and any $\phi \in H^1(\Gamma)$; \item $L$ has
eigenvectors but not a single one belongs to $H^s(\Gamma)$ for any
$s>0.$
\end{itemize}
\end{thm}

Later we will discuss examples of relaxation enhancing flows on
manifolds. One of our examples is derived from a construction
going back to Kolmogorov \cite{Kolmogorov}, and yields a smooth
flow with discrete spectrum and rough eigenfunctions. This example
is even more striking than the ones we discussed here since the
spectrum is discrete. However, the construction is more technical
and is postponed till Section~\ref{flowex}.

\section{The Fluid Flow Theorem}\label{flows}

In this section we discuss applications of the general criterion
to various situations involving diffusion in a fluid flow. First,
we are going to prove Theorem~\ref{fluid}. Most of the results we
need regarding the evolution generated by incompressible flows are
well-known and can be found, for example, in \cite{MP} in the
Euclidean space case. There are no essential changes in the more
general manifold setting.

\begin{proof}[Proof of Theorem \ref{fluid}]
It is well known that the Laplace-Beltrami operator $\Delta$ on a
compact smooth Riemannian manifold is self-adjoint, non-positive,
unbounded, and has a discrete spectrum (see e.g. \cite{Berger}).
Moreover, it is negative when considered on the invariant subspace
of mean zero $L^2$ functions. Henceforth, this will be our Hilbert
space: $H\equiv  L^2(M) \ominus {\bf 1}.$ Obviously it is sufficient
to prove Theorem \ref{fluid} for $\phi_0\in H$ (i.e., when
$\overline{\phi}=0$). The Lipschitz class divergence free vector
field $u$ generates a volume measure preserving transformation
$\Phi_t(x),$ defined by
\begin{equation}\label{flow1}
\frac{d}{dt}\Phi_t(x) = u(\Phi_t(x)), \,\,\,\Phi_0(x)=x
\end{equation}
(see, e.g. \cite{MP}). The existence and uniqueness of solutions to
the system \eqref{flow1} follows from the well-known theorems on
existence and uniqueness of solutions to first order systems of ODEs
involving Lipschitz class functions. With this transformation we can
associate a unitary evolution group $U^t$ in $L^2(M)$ where $U^tf(x)
= f(\Phi_{-t}(x)).$ It is easy to see that $H$ is an invariant
subspace for this group. The group $U^t$ corresponds to $e^{iLt}$ in
the abstract setting of Section~\ref{abs}. Since $\tfrac{d}{dt} (U^t
f) = - u\cdot \nabla (U^t f)$ for all $f\in H^1(M)$ (the usual
Sobolev space on $M$), we see that the
group's self-adjoint generator, $L,$ is defined by $L = i u \cdot
\nabla$ on functions from $H^1(M)$. It is clear that condition
\eqref{Lcon1} is satisfied, since $\|u \cdot \nabla f \| \leq
C\|f\|_1$ for all $f \in H^1.$ It remains to check that the
condition \eqref{Lcon2} is satisfied, that is, $\|e^{iLt}f\|_1 \leq
B(t) \|f \|_1$. Notice that if $u(x)$ is Lipschitz, so is
$\Phi_t(x)$ for any $t.$ This follows from the estimate
(in the local coordinates and for a sufficiently small time $t$)
\[
|\Phi_t(x) - \Phi_t(y)| \leq |x-y| + \int\limits_0^t
|u(\Phi_s(x))-u(\Phi_s(y)|\,ds.
\]
Applying Gronwall's lemma, we get
\[
|\Phi_t(x) - \Phi_t(y)| \le |x-y| e^{\|u\|_{\Lip}t}
\]
for any $x,y.$ Now by the well-known results on change of variables
in Sobolev functions (see e.g. \cite{Ziemer}) and by the fact that
$\Phi_t$ is measure preserving, we have that
\[
\|U^tf\|_1 \leq C \|\Phi_{t}\|_{\Lip}\|f\|_1.
\]
This is exactly \eqref{Lcon2}, and the application of
Theorem~\ref{general} finishes the proof.
\end{proof}

The criterion of Theorem~\ref{general} can be applied to boundary
value problems as well. For the sake of simplicity, consider a
bounded domain $\Omega \subset \Rm^d$ with a $C^2$ boundary $\partial
\Omega.$ Let $u \in \Lip(\Omega)$ be a Lipschitz incompressible flow
such that its normal component is zero on the boundary: $u(x) \cdot
\hat{n}(x) =0$ for $x \in \partial \Omega$, with $\hat n(x)$ the
outer normal at $x$. Let $\phi^A(x,t)$ be the solution of
\begin{equation}\label{Neumann}
\partial_t \phi^A(x,t) +A u \cdot \nabla \phi^A(x,t) -\Delta
\phi^A(x,t)=0, \,\,\,\phi^A(x,0)=\phi_0(x), \,\,\,\frac{\partial
\phi^A}{\partial n} =0\,\,\,{\rm if }\,\,\,x \in \partial \Omega,
\end{equation}
where the Neumann boundary condition is satisfied in the trace
sense for almost every $t>0$. The existence of solution to
\eqref{Neumann} can be proved similarly to Theorem~\ref{exun}.


\begin{thm}\label{Neumann1}
In the Neumann boundary conditions setting, the flow $u \in
\Lip(\Omega)$ is relaxation enhancing according to the
Definition~\ref{relen} if and only if the operator $u \cdot
\nabla$ has no eigenfunctions in $H^1(\Omega)$ other than the
constant function.
\end{thm}

\begin{proof}
The proof is essentially identical to that of Theorem~\ref{fluid}.
The Laplacian operator with Neumann boundary conditions restricted
to mean zero functions plays a role of the self-adjoint operator
$\Gamma.$ The condition $u \cdot \hat{n} =0$ ensures that the
vector field $u$ generates a measure preserving flow $\Phi_t(x)$
via \eqref{flow1}, and thus the corresponding evolution group is
unitary. The estimates necessary for Theorem~\ref{general} to
apply are verified in the same way as in the proof of
Theorem~\ref{fluid}.
\end{proof}

To treat other types of boundary conditions, such as Dirichlet, one
needs to modify the relaxation enhancement definition. This is due
to the fact that in this case the solution of \eqref{maineq} always
tends to zero, rather than to the average of the initial datum.

\begin{definition}\label{relen2}
Let $\phi^A(x,t)$ solve evolution equation \eqref{Neumann}, but
with Dirichlet or more general heat loss type boundary conditions
\begin{equation}\label{Dirichlet}
 \frac{\partial \phi^A}{\partial n}(x,t) +
\sigma(x)\phi^A(x,t) =0, \,\,\,x \in \partial \Omega,
\,\,\,\sigma(x) \in C(\partial \Omega), \,\,\,\sigma (x)> 0
\end{equation}
where $n$ is the outer normal to $\partial\Omega$. Then we call
the divergence free flow $u \in \Lip(\Omega)$ relaxation enhancing
if for every $\tau$ and $\delta$ there exists $A(\tau, \delta)$
such that for $A > A(\tau, \delta)$ and
$\|\phi_0\|_{L^2(\Omega)}=1$ we have
$\|\phi^A(x,\tau)\|_{L^2(\Omega)} < \delta.$
\end{definition}

\noindent \it Remarks. \rm 1. Note that $\sigma(x) =\infty$ is not
excluded and may lead to the Dirichlet boundary conditions on a
part of the boundary. \\
2. The more general definition encompassing both Definitions
\ref{relen} and \ref{relen2} would assume that the solution tends
to a certain limit and define relaxation enhancement in terms of
speed up in reaching this limit.

It is well known that the Laplace operator with boundary
conditions \eqref{Dirichlet} is self-adjoint on the domain of
$H^2(\Omega)$ functions satisfying \eqref{Dirichlet} in the trace
sense in $L^2(\partial \Omega).$ We denote this operator
$\Delta_\sigma.$ The corresponding $H^1_\sigma(\Omega)$ space is
the domain of the quadratic form of $\Delta_\sigma,$ consisting of
all functions $\phi \in H^1(\Omega)$ such that $\int_{\partial
\Omega} \sigma(x) |\phi(x)|^2 ds$ is finite. In the Dirichlet
boundary condition case, formally corresponding to $\sigma(x)
\equiv \infty,$ we obtain the standard space $H_0^1(\Omega).$ Then
we have

\begin{thm}\label{Dirichlet1}
In the case of the heat loss boundary condition \eqref{Dirichlet},
the flow $u \in \Lip(\Omega)$ satisfying $u \cdot \hat{n}=0$ on
the boundary is relaxation enhancing according to
Definition~\ref{relen2} if and only if the operator $u \cdot
\nabla$ has no eigenfunctions in $H^1_\sigma(\Omega).$
\end{thm}

\begin{proof}
In the case of heat loss boundary conditions, it is well-known
that the principal eigenvalue of $\Delta_\sigma$ is positive, so
we can set $\Gamma = -\Delta_\sigma.$ Our space $H$ is now equal to
$L^2(\Omega).$ The rest of the proof remains the same as in
Theorem~\ref{Neumann1}.
\end{proof}

We note that the case of Dirichlet boundary conditions has been
treated in \cite{BHN} in a more general setting $u \in
L^\infty(\Omega)$ and without the assumption $u \cdot \hat{n}=0.$
The methods of \cite{BHN} are completely different from ours, and
rely on the estimates on the principal eigenvalue of
$-\Delta+Au\cdot \nabla$ and positivity of the corresponding
eigenfunction. In particular, as described in the introduction,
these methods do not seem to be directly applicable to the study
of the enhanced relaxation in the case of a compact manifold without
boundary or Neumann boundary conditions, where the principal
eigenvalue is always zero. The results of \cite{BHN} show that in
the Dirichlet boundary condition case, the flow $u$ is relaxation
enhancing in the sense of Definition~\ref{relen2} if and only if
$u$ does not have a first integral in $H^1_0(\Omega).$ In other
words, if and only if the operator $u \cdot \nabla$ does not have
an $H^1_0(\Omega)$ eigenfunction corresponding to the eigenvalue
zero. The discrepancy between this result and
Theorem~\ref{Dirichlet1} may seem surprising, but in fact the
explanation is simple.

\begin{proposition}\label{simple11}
Let $u \in \Lip(\Omega).$ If $\phi \in H^1(\Omega)$
($H^1_\sigma(\Omega)$) is an eigenfunction of the operator $u \cdot
\nabla$ corresponding to the eigenvalue $i \lambda,$ then $|\phi|
\in H^1(\Omega)$ ($H^1_\sigma(\Omega)$) and it is the first integral
of $u$, that is, $u \cdot \nabla |\phi|=0.$
\end{proposition}

\begin{proof}
The fact that $|\phi| \in H^1$ follows from the well-known
properties of Sobolev functions (see e.g. \cite{Evans}). A direct
computation using $u \cdot \nabla \phi = i \lambda \phi$ then
verifies that $u \cdot \nabla |\phi|=0.$
\end{proof}

As a consequence, when $\sigma\not\equiv 0$, the condition of no
$H^1_\sigma$ eigenfunctions in the statement of
Theorem~\ref{Dirichlet1} can be replaced by
the condition of no first integrals in $H^1_\sigma.$ In the settings
of Theorems \ref{fluid} and
\ref{Neumann1}, the above argument still applies but does not
allow to change their statements. Indeed --- on one hand the
operator $u \cdot \nabla$ always has eigenvalue zero with an
eigenfunction that is smooth, namely a constant. Existence of this
first integral, however, tells us nothing about relaxation
enhancement. On the other hand, existence of mean zero $H^1$
eigenfunctions need not guarantee the existence of a mean zero
first integral, as can be seen in the following well-known
example.
\smallskip

\it Example. \rm Let $M\equiv T^d$ be the flat $d$-dimensional
torus with period one. Let $\alpha$ be a $d$ dimensional constant
vector generating irrational rotation on the torus (that is, we
assume that components of $\alpha$ are independent over the field
of rationals). It is well known that the flow generated by the
constant vector field $\alpha$ is ergodic but not weakly mixing.
The self-adjoint operator $L = i \alpha \cdot \nabla$ has
eigenvalues $2\pi \alpha \cdot k,$ where $k$ are all possible
vectors with integer components. The corresponding eigenfunctions
are $e^{-2\pi ik \cdot x},$ $x \in T^d.$ Their absolute value is
$1,$ which is a first integral of $\alpha,$ but there are no other
first integrals. In particular, every non-constant eigenfunction
of $L$ corresponds to a non-zero eigenvalue. Thus, this flow is
not relaxation enhancing even though it has no first integrals
other than a constant function.
\smallskip

Finally, we show that the $L^2$ norm in the
Definitions~\ref{relen}, \ref{relen2} can be replaced by other
$L^p$ norms with $1 \leq p \leq \infty$ without any change to the
statements of Theorems~\ref{fluid}, \ref{Neumann1},
\ref{Dirichlet1}. This result is important for applications to
quenching in reaction-diffusion equations.

\begin{thm}\label{anyp}
Theorems~\ref{fluid}, \ref{Neumann1}, \ref{Dirichlet1} remain true
if, in Definitions~\ref{relen}, \ref{relen2},
``$\|\phi_0\|_{L^2}=1$'' is replaced by ``$\|\phi_0\|_{L^p}=1$'' and
``$\|\phi^A(\tau)-\bar\phi\|_{L^2} < \delta$'' (resp.
``$\|\phi^A(\tau)\|_{L^2} < \delta$'') by
``$\|\phi^A(\tau)-\bar\phi\|_{L^q} < \delta$'' (resp.
``$\|\phi^A(\tau)\|_{L^q} < \delta$'') for any $p,q \in [1,\infty].$
\end{thm}

For the sake of consistency of notation, we will consider the
compact manifold case. The case of a domain $\Omega$ with
Dirichlet or heat-loss boundary conditions is handled similarly
(see below).

We start with the proof of a general $L^1\to L^\infty$ estimate
for solutions of
\begin{equation} \label{anyflow}
\psi_t + v \cdot \nabla \psi - \Delta \psi =0
\end{equation}
on a compact manifold $M$. The point is that this estimate will be
independent of the incompressible flow $v$ and so, in particular,
of the amplitude $A$ in \eqref{maineq}. It appeared, for example,
in \cite{FKR}, where the domain was a strip in $\Rm^2$. The
crucial ingredient of the proof was a Nash inequality. In the
general case, we follow a part of the argument, but our proof of
the corresponding inequality \eqref{nashineq} is different.

\begin{lemma} \label{nash}
For any smooth Riemannian manifold $M$ of dimension $d$ and any
$\eps\ge 0$ (resp. $\eps> 0$) if $d\ge 3$ (resp. $d=2$), there
exists $C=C(M,\eps)>0$ such that for any incompressible flow
$v\in\Lip(M)$ and any mean zero $\phi_0\in L^2(M)$, the solution
of \eqref{anyflow} and $\phi(x,0)=\phi_0(x)$ satisfies
\begin{equation} \label{evolbound}
\|\phi(x,t)\|_{L^\infty(M)} \le C t^{-d/2-\eps} \|\phi_0\|_{L^1(M)}.
\end{equation}
\end{lemma}

\begin{proof}
First note that by H\" older and Poincar\' e inequalities we have
for any mean zero $\psi\in H^1(M)$, any $p\ge (d+2)/4$ (if $d=2$,
then for any $p> (d+2)/4$), and some $C_p$,
\[
\|\psi\|_{L^2}^2 \le \|\psi\|_{L^1}^{1/p}
\|\psi\|_{L^{(2p-1)/(p-1)}}^{(2p-1)/p} \le C_p \|\psi\|_{L^1}^{1/p}
\|\nabla \psi\|_{L^2}^{(2p-1)/p}.
\]
That is,
\begin{equation} \label{nashineq}
\|\nabla \psi\|_{L^2}^{2}\ge C_q
\|\psi\|_{L^2}^{2+q}\|\psi\|_{L^1}^{-q}
\end{equation}
for $q\equiv 2/(2p-1)$ so that $q\le 4/d$ if $d\ge 3$ and $q<2$ if
$d=2$.

After multiplying \eqref{anyflow} by $\phi$ and integrating over $M$
we obtain for $t>0$
\begin{equation} \label{l2prime}
\frac d{dt}\|\phi\|_{L^2}^2 = -2\|\nabla\phi\|_{L^2}^2 \le -2C_q
\|\phi\|_{L^2}^{2+q}\|\phi\|_{L^1}^{-q} \le
-2C_q\|\phi\|_{L^2}^{2+q}\|\phi_0\|_{L^1}^{-q}.
\end{equation}
The last inequality follows from the positivity of and the
preservation of $L^1$ norms of solutions of \eqref{anyflow} with
initial conditions $\phi_{0,\pm}\equiv \max\{\pm\phi_0,0\}$, which
shows that $\|\phi\|_{L^1}$ is non-increasing.

Next we divide \eqref{l2prime} by $-\|\phi\|_{L^2}^{2+q}$ and
integrate in time to obtain $\|\phi(x,t)\|_{L^2}^{-q} \ge
qC_qt\|\phi_0\|_{L^1}^{-q}$. This in turn gives with a new $C_q$
\begin{equation} \label{evolbound2}
\|\phi(x,t)\|_{L^2} \le C_qt^{-1/q} \|\phi_0\|_{L^1}.
\end{equation}
Hence we have showed that $\|\calP_t(v)\|_{L^1\to L^2}\le
C_qt^{-1/q}$ where $\calP_t(v)$ is the solution operator for
\eqref{anyflow}. But since $\calP_t(v)=(\calP_t(-v))^*$ is the
adjoint of the operator $\calP_t(-v)$, which satisfies the same
bound, we obtain
\[
\|\calP_{2t}(v)\|_{L^1\to L^\infty} \le \|\calP_t(v)\|_{L^1\to L^2}
\|\calP_t(v)\|_{L^2\to L^\infty} = \|\calP_t(v)\|_{L^1\to L^2}
\|\calP_t(-v)\|_{L^1\to L^2} \le C_q^2t^{-2/q}
\]
which is \eqref{evolbound}.
\end{proof}

\begin{proof}[Proof of Theorem \ref{anyp}]
Assume for simplicity that the total volume of $M$ is equal to
one. Then it is clear that Lemma \ref{nash} also holds with
\begin{equation} \label{evolboundpq}
\|\phi(x,t)\|_{L^p(M)} \le C t^{-d/2-\eps} \|\phi_0\|_{L^q(M)}
\end{equation}
in place of \eqref{evolbound}, with any $p,q\in[1,\infty]$ and the
same $C$.

Assume now that we know that for some $u$ and $p,q\ge 1$, given
any $\tau,\delta>0$, we can find $A_{p,q}(\tau, \delta)$ such that
$\|\phi^A(x,\tau)\|_{L^p} < \delta$ for any $A> A_{p,q}(\tau,
\delta)$ and any mean zero $\phi_0\in L^2(M)$ with
$\|\phi_0\|_{L^q}=1$. Take any other $p',q' \geq 1.$ Then for any
$A> A_{p,q}(\tau, \delta)$,
\begin{equation*}
\|\phi^A(x,3\tau)\|_{L^{p'}}  \le C \tau^{-d/2-\eps}
\|\phi^A(x,2\tau)\|_{L^p}  \le \del C \tau^{-d/2-\eps}
\|\phi^A(x,\tau)\|_{L^q} \le \del(C \tau^{-d/2-\eps})^2
\|\phi_0\|_{L^{q'}}
\end{equation*}
This shows that when $A_{p,q}(\tau,\del)$ exists for some $p,q$
and all $\tau,\del$, for any $p',q',\tau,\del$ we have
$A_{p',q'}(\tau,\del) = A_{p,q}(\tau/3,\del
C^{-2}(\tau/3)^{d+2\eps})$ and so $A_{p',q'}(\tau,\del)$ exists
for all $\tau,\del$. That is, Definition \ref{relen} describes the
same class of flows regardless of which $L^p\to L^q$ decay it
addresses. This finishes the proof.
\end{proof}

We note that in the case of a bounded domain $\Omega$ with
Dirichlet boundary conditions, the proof is identical. When we
have heat-loss boundary conditions, the only change is that the
equality in \eqref{l2prime} reads
\[
\frac d{dt}\|\phi\|_{L^2}^2 = -2(\|\nabla\phi\|_{L^2}^2 +
\|\sigma^{1/2}\phi\|_{L^2(\partial\Omega)}^2)
\]
and the Poincar\' e inequality is replaced by
\[
\|\psi\|_{L^{(2p-1)/(p-1)}}\le C_p (\|\psi\|_{L^2} + \|\nabla
\psi\|_{L^2}) \le (C_p+\lambda_0^{-1/2}) (\|\nabla\psi\|_{L^2} +
\|\sigma^{1/2}\psi\|_{L^2(\partial\Omega)})
\]
which is due to the Sobolev inequality and the fact that the
principal eigenvalue $\lambda_0$ of the Laplacian on $\Omega$ with
heat-loss boundary conditions is positive.

\section{Examples of relaxation enhancing flows}\label{flowex}

Here we discuss examples of flows that are relaxation enhancing.
Most of the results in this section are not new and are provided for
illustration purposes. According to Theorem~\ref{fluid} a flow
$u\in \Lip(M)$ is relaxation enhancing if
all of its eigenfunctions are not in $H^1(M)$. One natural class
satisfying this condition is weakly mixing flows -- for which the
spectrum is purely continuous. Examples of weakly mixing flows on
$\Tm^2$ go back to von Neumann \cite{vN} and Kolmogorov
\cite{Kolmogorov}. The flow in von Neumann's example is
continuous; in the construction suggested by Kolmogorov the flow
is smooth. The technical details of the construction have been
carried out in \cite{Sklover}; see also \cite{Katok} for a review.
Recently, Fayad \cite{Fayad1} generalized this example to show
that weakly mixing flows are generic in a certain sense. For more
results on weakly mixing flows, see for example
\cite{Fayad2,Katok}. To describe the result in \cite{Fayad1} in
more detail, let us recall that a vector $\alpha$ in $\Rm^d$ is
called $\beta$-Diophatine if there exists a constant $C$ such that
for each $k\in\Zm^d \setminus \{0\}$ we have
\[
\inf_{p\in\Zm}|\langle \alpha, k \rangle +p|\ge
\frac{C}{|k|^{d+\beta}}.
\]
The vector $\alpha$ is Liouvillean if it is not Diophantine for
any $\beta>0.$ The Liouvillean numbers (and vectors) are the ones
which can be very well approximated by rationals.

\it Example 1. \rm Consider the flow on a torus $\Tm^{d+1}$ that
is a time change of a linear translation flow:
\begin{equation}\label{shift-time}
\frac{dx}{dt}=\frac{\alpha}{F(x,y)},~~\frac{dy}{dt}=\frac{1}{F(x,y)},
~~(x,y)\in\Tm^{d+1}
\end{equation}
with a smooth positive function $F(x,y)$. Such flows have a unique
invariant measure $d\mu=F(x,y)dx dy.$ Let us denote by
$C^r(\Tm^d,\Rm^+)$ the set of $C^r$ functions on the torus that are
positive. We have
\begin{proposition}[\cite{Fayad1}]\label{Fayad11}
Assume the irrational vector $\alpha \in \Rm^d$ is not
$\beta$-Diophantine, for some $\beta \in \Rm^+ \cup +\infty$.
Then, for a dense $G_\delta$ of functions $F$ in
$C^{\beta+d}(\Tm^d,\Rm^+)$ the flow \eqref{shift-time} is weakly
mixing (for the unique invariant measure $F(x,y)dx dy$).
\end{proposition}

To obtain examples of relaxation enhancing flows, we can now
consider the generic flows of Proposition~\ref{Fayad11} on the
torus with a metric such that the volume element is $F(x,y)dxdy.$
Alternatively, we can just view a problem in a weighted space and
consider the operator $\frac{1}{F} \nabla (F \nabla ),$ which is
self-adjoint on $L^2(\Tm^{d+1},Fdxdy),$ instead of $\Delta.$ It is
also straightforward to obtain more physical examples of the
relaxation enhancing flows that are incompressible with respect to
the usual flat metric. Indeed, assume for the sake of simplicity
that we are working with a unit torus and that the total integral
of $F$ is also one. Then it is not difficult to construct an
explicit measure preserving invertible transformation $Z$ from
$\Tm^{d+1}(F(x,y)dxdy)$ to $\Tm^{d+1}(dpdq),$ as smooth as the
function $F$ (for general results on existence of such maps see
\cite{Moser}). If we denote $w(x,y) = (\alpha/F(x,y),1/F(x,y))$
the vector field in \eqref{shift-time}, then the vector field
$u(p,q) = Z \circ w \circ Z^{-1}$ is going to be incompressible.
Moreover, the unitary evolutions generated by $w$ and $u$ in
$L^2(F(x,y)dxdy)$ and $L^2(dpdq)$ respectively are unitary
equivalent and so have the same spectra.

We now describe an example of a different class of flows to which
Theorem \ref{fluid} applies. Namely, we will sketch a construction
of a smooth incompressible flow $u(p,q)$, $\nabla\cdot u=0$, on a
torus $\Tm^2$ such that it has a purely discrete spectrum but none
of the eigenfunctions are in $H^1(\Tm^2)$. We could not find an
exact statement regarding the existence of such flows in the
literature, although the idea of the construction appears in
\cite{Kolmogorov} and the result is presumably well known in the
dynamical systems community. In particular, it follows in a fairly
direct way for example from considerations in \cite{Anosov,Katok}.
We briefly sketch the construction, without presenting well known
technical details.

\it Example 2. \rm Let us denote by $\Phi_t^u$ the flow on the
torus generated by $u$ and by $U^t$ the flow on $L^2(\Tm^2)$
generated by $\Phi_t^u$: $(U^t f)(x)=f(\Phi_{-t}^u(x))$.
\begin{proposition}\label{thm-flow-example}
There exists a smooth incompressible (with respect to the Lebesgue
measure) flow $u(x,y)$ on a two-dimensional torus $\Tm^2$ so that
the corresponding unitary evolution $U^t$ has a discrete spectrum
on $L^2(\Tm^2)$ but none of the eigenfunctions of $U$ are in
$H^1(\Tm^2)$.
\end{proposition}
\begin{proof}
The example will be given by a
flow of type \eqref{shift-time}, with $d=2$ and appropriately
chosen $\alpha$ and $F(x,y).$ We remark that while the form
\eqref{shift-time} may seem quite special, in fact any analytic
flow in two dimensions with an integral invariant can be mapped
analytically to the linear translation flow \eqref{shift-time}
with some $\alpha,$ $F$ (see e.g. \cite{Kolmogorov}). The idea of
the construction is to find a smooth flow \eqref{shift-time} which
can be mapped to a constant flow $(\alpha,1)$ by a measure
preserving map $S$ with very low regularity properties. Since the
eigenfunctions of the constant flow are explicitly computable, we
can compute the eigenfunctions of the original flow. Due to the
roughness of $S,$ these will prove highly irregular. To obtain an
incompressible flow, we will then proceed as in the first example.

In order to find such a flow $w$, we start with a smooth periodic
function $Q\in C^\infty(\Sm^1)$ and an irrational number
$\alpha\in\Rm$ so that the homology equation
\begin{equation}\label{homology}
R(\xi+\alpha)-R(\xi)=Q(\xi)-1,~~~\xi\in{\mathbb S}^1,
\end{equation}
has a solution $R(\xi)$ that is very rough.  Note that for
(\ref{homology}) to have a measurable solution the function
$Q(\xi)$ should satisfy the normalization \cite{Anosov}
\[
\int_0^1 Q(\xi)d\xi=1.
\]
The following Proposition is a particular case of Theorem 4.5 of
\cite{Katok}.
\begin{proposition}\label{prop-meas}
Let $\alpha$ be a Liouvillean irrational number. There exists a
$C^\infty(\Sm^1)$ function $Q(\xi)$ so that the homology equation
(\ref{homology}) has a unique (up to an additive constant)
measurable solution $R(\xi):{\Sm}^1\to\Rm$ such that for any
$\lambda\in\Rm_*=\Rm\backslash\{0\}$, the function
$R_\lambda(\xi)=e^{i\lambda R(\xi)} $ is discontinuous everywhere.
\end{proposition}
Note that without loss of generality we may assume that $Q(\xi)$
is positive -- otherwise we choose $M$ so that $Q(\xi)+M>1$ and
consider a rescaled function $Q_M(\xi)=(M+Q(\xi))/(M+1)$. Then the
function $R_M(\xi)=R(\xi)/(M+1)$ is the solution of
(\ref{homology}) with $Q_M$ on the right side and, of course,
$R_M(\xi)$ has the same properties as $R(\xi)$.

Given a Liouvillean irrational number $\alpha$ and a function
$Q(\xi)$ that satisfies the conclusion of Proposition
\ref{prop-meas} we define a function $F(x,y)$ on the torus $\Tm^2$
as follows. Choose $m>0$ so that $m<\min Q(s)$ and a smooth
function $\psi(y)\ge 0$ such that
\begin{equation}\label{psi-int1}
\int_0^1\psi(y)dy=1.
\end{equation}
and, in addition,
\begin{equation}\label{psi-equal}
\hbox{$\psi(y)=0$ for $0\le y\le y_0$ and $y_1\le y\le 1$ with
$y_0$ close to zero and $y_1$ close to one.}
\end{equation}
The choice of $m$ ensures that the function
\begin{equation}\label{def-G}
F(x,y)=m+\psi(y)(Q(x-\alpha y)-m),~~0\le x,y\le 1
\end{equation}
is positive -- then we extend $F(x,y)$ periodically in both
variables to the whole plane $\Rm^2$. The resulting function is
smooth because of (\ref{psi-equal}) and, in addition, it has total
mass equal to one. The normalization (\ref{psi-int1}) implies that
the functions $F$ and $Q$ are related by
\begin{equation}\label{Q-def}
Q(\xi)=\int_0^1{F(\xi +\alpha z,z)}dz.
\end{equation}
Now, the required transformation $S:(x,y)\to (X,Y)$ is defined by
\cite{Kolmogorov,Sternberg}
\begin{equation}\label{transform}
X(x,y)=x+\alpha(Y(x,y)-y),~~~~Y(x,y)=T(x-\alpha y,y)+R(x-\alpha y)
\end{equation}
with the function $R(x)$ that satisfies (\ref{homology}), and
$T(x,y)$ defined by
\[
T(x,y)=\int_0^y F(x+\alpha z,z){dz}.
\]
Note that the transformation \eqref{transform} implies that $x -
\alpha y = X -\alpha Y,$ and so it preserves the flow
trajectories. The homology equation (\ref{homology}) together with
the definition (\ref{def-G})  of the function $F(x,y)$ imply that
$S$ is well-defined as a mapping $\Tm^2\to \Tm^2.$ It is also
straightforward to check that $S$ maps the flow $w$ onto the
uniform flow $w_{unif}=(\alpha,1)$. One can also verify that $S$
is invertible with measurable inverse, and is measure preserving:
\begin{equation}\label{S*}
\int[S^*f](x,y)F(x,y)dxdy=\int f(S(x,y))F(x,y)dxdy=\int f(X,Y)dXdY
\end{equation}
for any function $f\in C(\Tm^2).$ Hence, $S^*$ may be extended as
an operator $L^2(dxdy)\to L^2(d\mu)$ with the preservation of the
corresponding norms. We conclude that the unitary evolutions
$U_w^t$ and $U^t_{unif}$ generated by the flow $w$ given by
\eqref{shift-time} and the uniform flow $w_{unif}$, respectively,
are conjugated by means of a unitary transformation
$S^*:~L^2(\Tm^2,dXdY)\to L^2(\Tm^2,d\mu)$ -- we have
$U^t_{unif}=[S^*]^{-1}U^t_{w} S^*$. It follows that $U^t_w$ and
$U^t_{unif}$ have the same spectrum: $\lambda_{nl}=2\pi in\alpha
+2\pi i l, ~~~l,n\in{\mathbb Z}.$ It also follows that the
eigenfunctions of the operator $U_w$ may be written as
\begin{eqnarray}\label{eigen-u2}
&&\psi_{nl}^w(x,y)=e^{2\pi in X(x,y)+2\pi i l Y(x,y)}=e^{2\pi
in(x-\alpha y+\alpha Y(x,y))+2\pi i l Y(x,y)}
\\
&&\,\,\,\,\,\,\,\,\,\,\,\,\,\,\,\,\,\,\,\,\,\, =e^{2\pi in
(x-\alpha y)}e^{(2\pi i n\alpha +2\pi il)(T(x-\alpha
y,y)+R(x-\alpha y))} =\zeta(x,y)e^{(2\pi i n\alpha +2\pi
il)R(x-\alpha y)}\nonumber
\end{eqnarray}
with a smooth function $\zeta(x,y)\in C^\infty([0,1]^2)$ (note
that the function $\zeta(x,y)$ is not periodic in $(x,y)$). In
order to verify that $\psi_{nl}^w$ are not in $H^1(\Tm^2)$ it
suffices to check that the function
\[
\Theta_\lambda(x,y)=e^{i\lambda R(x-\alpha y)}= R_\lambda(x-\alpha
y)
\]
is not in $H^1([0,1]^2)$ for any real $\lambda\neq 0$. The
function $R_\lambda(s)$ is defined in the Proposition
\ref{prop-meas} and is everywhere discontinuous. If
$\Theta_\lambda(x,y)$ were in $H^1([0,1]^2)$, it would force
$R_\lambda(s)$ to be in $H^1(\Sm^1)$ and hence continuous but this
function is discontinuous everywhere. Therefore, the
eigenfunctions $\psi_{nl}^w$ cannot be in $H^1(\Tm^2).$

Finally, to obtain an incompressible flow, we introduce a smooth
transformation $Z:~(x,y)\to(p,q)$ be setting
\[
p=\int_0^x\bar {F}(s)ds,~~q=\frac{1}{\bar F(x)}\int_0^y
\!F(x,z)dz, \hbox{ where  } \bar F(x)=\int_0^1\!F(x,z)dz.
\]
It is immediate to verify that $Z$ maps the measure $d\mu$ onto
the Lebesgue measure $dpdq$. Hence, the evolution group generated
by the image $u(p,q)$ of the flow $w(x,y)$ will have the same
discrete spectrum as $U_w$. In addition, the eigenfunctions
$\psi_{nl}^w$ of $U_w$ are the images of the eigenfunctions
$\psi_{nl}^u$ of $u$ under $Z^*$:
$\psi_{nl}^w=Z^*\psi_{nl}^u=\psi_{nl}^u\circ Z$. As the functions
$\psi_{nl}^w$ are not in $H^1(\Tm^2)$ and the map $Z$ is smooth,
it follows that all the eigenfunctions of the incompressible flow
$u(p,q)$ are not in $H^1(\Tm^2)$. This finishes the proof of
Proposition \ref{thm-flow-example}.
\end{proof}

\section{Quenching in Reaction-Diffusion Equations}\label{quenching}

In this section we describe the application of our results to
questions of quenching in reaction-diffusion-advection
equations. We will consider the problem
\begin{equation}\label{readiff}
T^A_t(x,t) +A u \cdot \nabla T^A(x,t) - \Delta T^A(x,t)=f(T^A(x,t)),
\,\,\,\,\,T^A(x,0)=T_0(x)
\end{equation}
on a smooth compact Riemannian manifold $M$ with $T_0(x)\in[0,1]$.
Here $T$ is the (normalized) temperature of a premixed flammable
gas that is advected by the incompressible flow $Au(x)\in\Lip(M)$.
The nonlinear right hand side term accounts for temperature
increase due to burning and will be assumed to be of {\it
ignition} type. That is,
\begin{equation} \label{ignition}
\begin{split}
& \text{(i) $f(0)=f(1)=0$ and $f(T)$ is Lipschitz continuous on
$[0,1]$},
\\ & \text{(ii) $\exists\theta_0\in(0,1)$ such that $f(T)=0$ for
$T\in[0,\theta_0]$ and $f(T)> 0$ for $T \in (\theta_0,1)$}.
\end{split}
\end{equation}
This shows, in particular, that $T$ remains in $[0,1]$. The main
question will be under what conditions on $u$ one can always choose
$A$ large enough so that for some time $\tau>0$ we have
$\|T^A(x,\tau)\|_{L^\infty(M)}\le \tht_0$, that is, quenching
--- extinction of flames --- happens. Of course this question is
meaningless for certain initial data $T_0$. Namely, if
$\|T_0\|_{L^1}>\tht_0 {\rm vol}(M)$ or $\|T_0\|_{L^1}=\tht_0 {\rm
vol}(M)$ but $T_0\not\equiv\tht_0$, then it is easy to show using
\begin{equation} \label{timederiv}
\frac d{dt}\|T^A\|_{L^1} = \int f(T^A(x,t)) dx
\end{equation}
that $\|T^A\|_{L^1}$ must be strictly increasing with the limit
equal to the volume of $M$. This motivates the following definition.

\begin{definition} \label{quenchdef}
We say that $u$ is strongly quenching if for any nonlinearity $f$
as in \eqref{ignition}, and any solution $T^A$ of \eqref{readiff}
with initial datum $T_0(x)\in[0,1]$ with $\|T_0\|_{L^1(M)}<\tht_0
{\rm vol}(M)$, there exists $A(T_0,f)$ such that if $A>A(T_0,f)$,
then for some $\tau>0$ one has $\|T^A(x,\tau)\|_{L^\infty(M)}\le
\tht_0$.
\end{definition}

Then we have

\begin{theorem} \label{quenchthm}
An incompressible flow $u\in\Lip(M)$ is strongly quenching if and
only if it is relaxation enhancing.
\end{theorem}

\begin{proof}
Assume that the volume of $M$ is one. First, $u$ is strongly
quenching when it is relaxation enhancing. Indeed, assume $u$ is
relaxation enhancing and let $c$ be the Lipshitz constant for $f$
so that $f(T)\le cT$. If $\phi^A$ solves \eqref{maineq} with
$\phi_0\equiv T_0$, then $T^A(x,t)\le e^{ct}\phi^A(x,t)$ by the
comparison principle. But this means
\[
\|T^A(x,\tau)\|_{L^\infty}\le e^{c\tau}
\|\phi^A(x,\tau)\|_{L^\infty} \le e^{c\tau} (\bar
T_0+\|\phi^A(x,\tau)-\bar\phi\|_{L^\infty})
\]
which can be made as close to $\bar T_0=\|T_0\|_{L^1(M)}<\theta_0$
as we wish by taking small enough $\tau,\del>0$ and
$A>A_{1,\infty}(\tau,\del)$ from the proof of Theorem \ref{anyp}.
Since $f$ was arbitrary, it follows that $u$ is strongly
quenching.

Hence, we are left with proving that $u$ being  strongly quenching
implies that $u$ is relaxation enhancing. Assume this is not the
case, that is, there exists an ignition nonlinearity $f$, a mean
zero $\phi_0\in L^2(M)$ and $\tau,\del\in(0,1)$ such that for all
$A<\infty$ and $t\le\tau$, the solution $\phi^A$ of \eqref{maineq}
satisfies $\|\phi^A(x,t)\|_{L^1}>\del$. We can assume without loss
of generality that $\|\phi_0\|_{L^\infty}\le 1$.

Let $T_0\equiv \tht_1+\gamma\phi_0$ where $\gamma\equiv \tfrac
13\min\{\tht_0,1-\tht_0 \}$ and $\tht_1\in(\tht_0-\beta,\tht_0)$
with $\beta\equiv\tfrac 18 \min\{\gamma\del, \tau \del\kappa
,\tht_0\}$ and $\kappa\equiv\min\{f(T)\,|\, T\in
[\tht_0+\tfrac{\gamma\del}8,\tfrac13 (2+\tht_0)]\}>0$. Note that
$T_0(x)\in\tfrac 13[\tht_0,1+2\tht_0]$.

For each $t\le\tau$ let $B_t\equiv\{x\,|\,
\phi^A(x,t)\ge\tfrac\del 4\}$. Since $\phi^A(x,t)$ is mean zero
with $L^1$ norm more than $\del$, and $|\phi^A(x,t)| \leq 1,$ we
must have $|B_t|\ge\tfrac \del 4$ for each $t\le\tau$ (recall that
$M$ has volume one). Thus
\[
T^A(x,t)\ge\tht_1+\tfrac{\gamma\del} 4\ge
\tht_0+\tfrac{\gamma\del} 8
\]
for $x\in B_t$. If for some $t\le\tau$ there is a set
$B'_t\subseteq B_t$ with $|B'_t|\ge\tfrac \del 8$ and $T^A(x,t)\ge
\tfrac 13(2+\tht_0)$ for $x\in B'_t$, then since
$\psi^A\equiv\tht_1+\gamma\phi^A\le T^A$ by the comparison
principle (because $\psi^A$ satisfies \eqref{maineq} with
$\psi^A(x,0)=T_0(x)$) and $0<\inf\{T_0(x)\}\le\psi^A\le
\|T_0\|_{L^\infty}\le \tfrac 13(1+2\tht_0)$ by the maximum
principle,
\[
\|T^A(x,t)\|_{L^1}\ge \tfrac \del 8 [\tfrac 13(2+\tht_0)-\tfrac
13(1+2\tht_0)] + \|\psi^A(x,t)\|_{L^1}  \ge \tfrac{\gamma\del}8
+\tht_1
>\tht_0.
\]
This is a contradiction because then $u$ cannot be strongly quenching
by the argument before Definition \ref{quenchdef}.

Therefore for each $t\le\tau$ there must be a set $B''_t\subseteq
B_t$ such that $|B''_t|\ge\tfrac \del 8$ and $T^A(x,t)\in
[\tht_0+\tfrac{\gamma\del}8, \tfrac 13(2+\tht_0)]$ for $x\in
B''_t$. But then $f(T^A(x,t))\ge\kappa$ for $x\in B''_t$, and so
\eqref{timederiv} gives
\[
\|T^A(x,\tau)\|_{L^1} \ge \tht_1 + \tau\tfrac{\del\kappa}8>\tht_0
\]
and we have a contradiction again. Hence $u$ has to be relaxation
enhancing.
\end{proof}

Applications of our results on relaxation enhancement to quenching
on infinite domains (where front propagation can occur as well)
will be considered elsewhere.

\smallskip
\noindent {\bf Acknowledgement.} PC has been partially supported
by the NSF-DMS grant 0202531. AK and LR have been supported in
part by the Alfred P. Sloan Research Fellowships. In addition, AK
and AZ have been supported by the NSF-DMS grant 0314129. AK and AZ
thank the University of Chicago for its hospitality in the spring
of 2005. The authors thank Mark Freidlin, Alex Furman, Stas
Molchanov and Laurent Saloff-Coste for interesting discussions and
for pointing out some relevant literature.

\end{document}